\documentclass[11pt, english]{smfart}
\def\today{16.10.06} 
\usepackage{amsmath,amsfonts,amsthm,amssymb,amscd}
\usepackage[francais,english]{babel}
\newcommand{\be}{\begin{equation}}
\newcommand{\ee}{\end{equation}}
\newcommand{\bi}{\begin{itemize}}
\newcommand{\ei}{\end{itemize}}

\theoremstyle{plain} \newtheorem{theorem}{Theorem}[section]
\newtheorem{lemma}[theorem]{Lemma}
\newtheorem{proposition}[theorem]{Proposition}
\newtheorem{corollary}[theorem]{Corollary} 
\newtheorem{definition}[theorem]{Definition} \theoremstyle{remark}
\newtheorem{remark}[theorem]{Remark}
\newtheorem{example}[theorem]{Example}
\newtheorem{exercise}[theorem]{Exercise} 
\newcommand{\R}{{\mathbb R}} \newcommand{\U}{{\mathcal U}}

\newcommand{\Z}{{\mathbb Z}}
\newcommand{\N}{{\mathbb N}}
\newcommand{\e}{{\varepsilon}}

\newcommand{\Zb}{{\bar \Z}} \newcommand{\Nb}{{\bar \N}}

\newcommand{\V}{{\mathcal V}}
\newcommand{\W}{{\mathcal W}}
\newcommand{\Ps}{{\mathcal P_{s}}}

\newcommand{\ls}{{\mathcal L_{s}}}
\newcommand{\Hs}{{\mathcal H^s}}
\newcommand{\Qs}{{\mathcal Q_s}}

\def\im{{\rm i}}

\newcommand{\C}{\mathbb{C}}
\newcommand{\T}{\mathcal{T}}
\newcommand{\Tt}{\mathbb{T}}

\newcommand{\St}{\tilde S}

\newcommand{\va}[1]{|#1|}
  \newcommand{\vak}{|k|}
  \newcommand{\vaj}{|j|}
  \newcommand{\val}{|l|}
  
  \newcommand{\vat}{|t|}
  \newcommand{\Va}[1]{\left|#1\right|}

\def\norma#1{\left\| #1\right\|}

\def\sleq{\leq\kern-6pt \cdot\null\hskip4pt}

\numberwithin{equation}{section}

\begin{document}

\author{ Beno\^it Gr\'ebert} \title{Birkhoff Normal Form and Hamiltonian PDEs}
\alttitle{Forme normale de Birkhoff et EDP Hamiltoni\`ennes}
\date{\today}

\begin{abstract}
    These notes are based on lectures held at the Lanzhou university 
    (China) during a CIMPA summer school in july 2004 but benefit 
    from recent devellopements. 
    Our aim is to explain some normal form technics 
    that allow to study the long time behaviour of the solutions of 
    Hamiltonian perturbations of integrable systems. We are in particular 
    interested with stability results. 
    
    Our approach is centered on the Birkhoff normal form theorem 
    that we first proved in finite dimension. Then, 
    after giving some exemples of Hamiltonian PDEs, we present an 
    abstract 
    Birkhoff normal form theorem in infinite dimension and discuss the dynamical 
    consequences for Hamiltonian PDEs.
\end{abstract}

\begin{altabstract}
    Ces notes sont bas\'ees sur un cours donn\'e \`a 
    l'universit\'e de Lanzhou (Chine) durant le mois de juillet 2004 
    dans le cadre d'une \'ecole d'\'et\'e organis\'ee par le CIMPA. 
    Cette r\'edaction
     b\'en\'eficie aussi de d\'eveloppements plus r\'ecents. Le 
    but est d'expliquer certaines techniques de forme normale qui 
    permettent d'\'etudier le comportement pour des temps longs des 
    solutions de perturbations Hamiltoni\`ennes de syst\`emes 
    int\'egrables. Nous sommes en particulier int\'eress\'es par des 
    r\'esultats de stabilit\'e.
    
    Notre approche est centr\'ee sur le th\'eor\`eme de forme normale 
    de Birkhoff que nous rappelons et d\'emontrons d'abord en 
    dimension finie. Ensuite, apr\`es avoir donn\'e quelques exemples 
    d'EDP Hamiltoni\`ennes, nous d\'emontrons un 
    th\'eor\`eme de forme normale 
	de Birkhoff en dimension infinie et nous en discutons les 
	applications \`a la dynamique des EDP Hamiltoni\`ennes.

\end{altabstract}    
\frontmatter
\maketitle
\newpage
\mainmatter
\tableofcontents

\setcounter{section}{0}

\section{Introduction}
The class of Hamiltonian systems close to integrable system contain 
most of the important physic models.
Typically a Hamiltonian system in finite dimension  reads (cf. section 
\ref{formalism}) 
$$
\left\{ \begin{array}{c}
\dot q_{j}=\frac{\partial H}{\partial p_{j}}, \quad j=1,\ldots,n  \\
\dot p_{j}=-\frac{\partial H}{\partial q_{j}}, \quad j=1,\ldots,n 
\end{array}\right.
$$
where the Hamiltonian $H$ is a smooth fonction from $\R^{2N} $ to 
$\R$. In these lectures we are interesred in the case where $H$ 
decomposes in $H=H_{0}+\epsilon P$, $H_{0}$ being integrable in the 
sense that we can "integrate" the Hamiltonian system associated to 
$H_{0}$ (cf. section \ref{int}), $P$ being the perturbation and 
$\epsilon$ a small parameter. This framework contains a lot of 
important example of the classical mechanics. If we allow the number 
of degree 
of freedom, $N$, to grow to 
infinity, then we arrive in the world of quantum mechanics and the 
corresponding equations are typically nonlinear partial differential equations 
(PDEs). Again a lot of classical examples are included in this 
framework like, for instance, the nonlinear wave equation, the 
nonlinear Schr\"odinger equation or the Korteweg-de Vries equation (cf. 
section \ref{ex}). 

The historical example (in finite dimension) is given by the celestial 
mechanics: More than 300 years ago Newton gaves the evolution 
equation for a system of $N$ heavy bodies under the action of the 
gravity.\\ When $N=2$, Kepler gaves the solution, the bodies describe 
ellipses. Actually for $N=2$ 
the system is integrable.\\ As soon as $N\geq 3$ the system leaves the 
integrable world and we do not know the expression of the general 
solution. Nevertheless if we consider the celestial system composed by 
the Sun (S), the Earth (E) and Jupiter (J) and if we neglect the 
interaction between J and E, then the system is again integrable and 
we find quasiperiodic solution. Mathematically the solutions read 
$t\mapsto g(\omega_{1}t,\omega_{2}t,\omega_{3}t)$ where $g$ is a 
regular function from the torus $T^3=S^1\times S^1\times S^1$ to 
$\R^{18}$ (three postions and three moments in $\R^3$) and 
$\omega_{j}$, $j=1,2,3$ are frequencies. Visually J and E turn around 
S which turns around the center of mass. Notice that the trajectory 
(or orbit) is contained in the torus $g(\T^3)$ of dimension 3 and that 
this torus is invariant under the flow. On the other hand, 
if $(\omega_{1},\omega_{2},\omega_{3})$ are rationnaly independent, 
then the trajectory densely fills this torus while, if for instance 
the three frequencies are rationnally proportional, then the 
trajectoty is periodic and  
describes a circle included in  $g(\T^3)$.\\
Now the exact system S-E-J is described by a Hamiltonian 
$H=H_{0}+\epsilon P$ in which $H_{0}$ is the integrable Hamiltonian where 
we neglect the interaction E-J, $P$ takes into account this 
interaction and $\epsilon= \frac{\mbox{jupiter's mass + earth's 
mass}}{\mbox{sun's mass}}$ plays the rule of the small parameter.\\
Some natural questions arrive:
\begin{itemize}
\item Do invariant tori persist after this small perturbation?
\item At least are we able to insure stability in the sense that the 
planets remain in a bouded domain?
\item Even if we are unable to answer these questions for eternity, 
can we do it for very large -but finite- times?
\end{itemize}
These questions have interested a lot of famous mathematicians and 
physicists. In the 19-th century one tried to expand the solutions in 
perturbative series: $u(t)=u_{0}(t)+\epsilon u_{1}(t)+\epsilon^2 
u_{2}(t)+\ldots$, the term $u_{k+1}$ being determined by an equation 
involving $u_{0},\ldots,u_{k}$. Unfortunatly this series does not 
converge. This convergence problems seemed so involved that, at the 
principle of the 20-th century, most of scientist believed in the 
\textit{ergodic hypothesis}: typically, after arbitrarily small perturbation, 
all the trajectories fill all the phase space and the stable 
trajectories are exceptionnal. Actually, H. Poinar\'e proved that 
a dense set of invariant tori are destroyed by an arbitrarily small perturbation.
Nevertheless, a set can be dense but very small and in 1954 A. N. 
Kolmogorov \cite{Kol} announced that the majority (in the measure sense) of tori 
survive (see section 7). The proof of this result was completed by V. 
Arnold \cite{A63} and J. 
Moser \cite{Mos62} giving birth to the KAM theory.  

In order to illustrate this result we can apply it to a simplified 
S-E-J system: we assume that the S-E-J system reduces to a Hamiltonian 
system with 3 degrees of freedom without symmetries (the symmetries of 
the true system complicates the pictures and generates 
degenerancies). In this  case, the KAM theorem says, roughly 
speaking (see theorem \ref{kam-fini} for a precise statement), that 
if $(\omega_{1},\omega_{2},\omega_{3})\in \mathcal C$, a Cantor set 
of $\R^3$ having a positive measure, or equivalently if the initial 
positions and moments are in a Cantor set, then the trajectory is quasi 
periodic. Since a Cantor set has an empty interior the condition 
$(\omega_{1},\omega_{2},\omega_{3})\in \mathcal C$ is not physical (no 
measurement could decide if this condition is verified or not). \\
The present
lectures will be centered on the Birkhoff normal form approach which 
does not control the solution 
for any times but does not require an undecidable hypothesis. 
In the case of our simplified S-E-J system, the Birkhoff normal form theorem 
says, roughly 
speaking, that 
having fixed an integer $M\geq 1$, and  $\epsilon <\epsilon_{0}(M)$ 
small enough,
to any initial datum corresponding to not rationaly dependent 
frequencies $(\omega_{1},\omega_{2},\omega_{3})$,
we can associate a torus such that the solution remains
 $\epsilon$-close to that torus during a lapse of time greater than 
 $1/\epsilon^M$ (see section 3 for a precise statement). Note that this result 
 can be physically sufficient if $1/\epsilon^M$ is greater than the 
 age of the universe.

The rational independence of the frequencies (one also says the 
nonresonancy) is of course essential in all this kind of perturbative 
theorems. Again we can illustrate this fact with our system S-E-J: 
suppose that, when considering the system without E-J interaction,
the three bodies are periodically align, the Earth being between 
Jupiter
and the Sun (notice that this implies that the frequencies 
$(\omega_{1},\omega_{2},\omega_{3})$ are rationaly dependent). When we 
turn on the interaction E-J, 
Jupiter will attract the Earth outside of 
its orbit periodically (i.e. when the three bodies are align or 
almost align), these accumulate small effects will force the earth 
to escape its orbit and thus the invariant torus will be destroyed.

\medskip

The generalisation of these results to the infinite diemensional case 
is of course not easy but it worth trying: The expected results may 
apply to nonlinear PDEs when they can be viewed as an
infinite dimensional Hamiltonian system (cf. section 5) and concern 
the long time behaviour of the solution, a very difficult and 
competitive domain. 

For a general overview on 
Hamiltonian PDEs, the reader may consult the recent monographies
by Craig \cite{Cr00}, by Kuksin \cite{K2}, by Bourgain 
\cite{Bo05} and by Kappleler 
and P\"oschel \cite{KaP}. In the present lectures we mainly focus on the 
extension of the
Birkhoff normal form theorem. Such extension was first (partially) achieved by Bourgain 
\cite{Bo96} and then by Bambusi \cite{Bam03}. The results stated in 
this text was first proved by Bambusi and myself in \cite{BG04}. The 
proof presented here and some generalisations benefit of a recent 
collaboration with Delort and Szeftel \cite{BDGS}.

After this general presentation, I give a brief outline 
of the next sections:
\begin{itemize}
\item \textbf{Section 2 }: Hamiltonian formalism in finite dimension.

We recall briefly the classical Hamiltonian formalism including: 
integrals of the motion, Lie transforations,
Integrability in the Liouville sense, action angle variables, 
Arnold-Liouville theorem (see for instance \cite{A1}  for a complete 
presentation).
\item \textbf{ Section 3} : The Birkhoff normal form theorem in finite dimension. 

We state and prove the Birkhoff normal form theorem and then present its 
dynamical consequences. Theses results are well known and the reader 
may consult \cite{MoS, HZ, KaP} for more details and generalizations.

\item \textbf{Section 4 }: A Birkhoff normal form theorem in infinite dimension.

We state a Birkhoff normal form theorem in infinite dimension and explain its 
dynamical consequences. In particular, results on the longtime
behaviour of the solutions  are discussed. This is the 
most important part of this course. A slightly more general abstract 
Birkhoff theorem in infinite dimension was obtained in \cite{BG04} 
and the dynamical consequences was also 
obtained there.

\item \textbf{Section 5} : Application to Hamiltonian PDEs.

Two examples of Hamiltonian PDEs are given: the nonlinear wave 
equation and the nonlinear Schr\"odinger equation. We then verify that 
our Birkhoff theorem and its dynamical consequences apply  to both 
examples.

\item \textbf{Section 6} : Proof of the Birkhoff normal form theorem in infinite 
dimension.

Instead of giving the proof of \cite{BG04},
we present a simpler proof using a class of 
polynomials first introduced in \cite{DS1}, \cite{DS2}. Actually we 
freely used notations and parts of proofs of theses three references.

\item \textbf{Section 7} : Generalisations and comparison with KAM 
type results.

In a first part we comment on some generalisations of our result.
In the second subsection, we try to give to the reader an idea on the KAM 
theory in both finite and infinite dimension.  
Then we  compare the Birkhoff approach with the KAM approach.
\end{itemize}    

{\it Acknowledgements: it is a great pleasure  to
thank D. Bambusi and J. M. Delort for helpful discussions on these 
notes.}

\section{Hamiltonian formalism in finite dimension}
\label{formalism} 
\subsection{Basic definitions}\label{basic}
We only consider the case where the \textbf{phase space} (or 
configuration space) is an open set, $M$, of $\R^{2n}$.
We denote by $J$ the canonical \textbf{Poisson matrix}, i.e.
$$J=
\left( \begin{array}{cc}
    0 & I_{n} \\
    -I_{n} & 0
\end{array}\right)\ .
$$
More generally, $J$ could be an antisymmetric matrix on $\R^{2n}$. 
All the 
theory can be extended to the case where the phase space is a 
$2n$ dimensional symplectic manifold.\\
A \textbf{Hamiltonian fonction}, $H$, is a regular real valued function on the phase 
space, i.e. $H\in C^{\infty}(M,\R)$. To $H$ we associate the 
\textbf{Hamiltonian vector field}
$$X_{H}(q,p)=J\nabla_{q,p}H(q,p)$$
where $\nabla_{p,q}H$ denotes the gradient of $H$ with respect to 
$p,q$, i.e.
$$\nabla_{q,p}H=
\left(\begin{array}{c}
    \frac{\partial H}{\partial q_{1}}\\
  \vdots  \\
  \frac{\partial H}{\partial q_{n}}  \\
  \frac{\partial H}{\partial p_{1}}  \\
    \vdots  \\
    \frac{\partial H}{\partial p_{n}} \\
\end{array}\right)\ , \quad 
X_{H}=
\left(\begin{array}{c}
    \frac{\partial H}{\partial p_{1}}\\
  \vdots  \\
  \frac{\partial H}{\partial p_{n}}  \\
  -\frac{\partial H}{\partial q_{1}}  \\
    \vdots  \\
    -\frac{\partial H}{\partial q_{n}} \\
\end{array}\right)\ .
$$
The associated \textbf{Hamiltonian system} then reads
$$
\frac{d}{dt}
\left(\begin{array}{c}
    q  \\
    p
\end{array}\right)=X_{H}(q,p)$$
or equivalently
$$
\left\{ \begin{array}{cccc}
\dot q_{j}&=&\frac{\partial H}{\partial p_{j}}, &\quad 
j=1,\ldots,n\, ,  \\
\dot p_{j}&=&-\frac{\partial H}{\partial q_{j}}, &\quad j=1,\ldots,n\, . 
\end{array}\right.
$$
The \textbf{Poisson bracket} of two Hamiltonian functions $F,G$ is a 
new Hamiltonian function $\{ F,G \}$ given  by
$$
\{ F,G \}(q,p)=\sum_{j=1}^n \ \frac{\partial F}{\partial q_{j}}(q,p)
\frac{\partial G}{\partial p_{j}}(q,p) -
\frac{\partial F}{\partial p_{j}}(q,p)\frac{\partial G}{\partial q_{j}}(q,p)
\ .$$
\subsection{A fundamental example: the harmonic oscillator} 
\label{harmonic}
Let $M=\R^{2n}$ and
$$H(q,p)=\sum_{j=1}^n \omega_{j} \frac{p_{j}^2+q_{j}^2}{2}$$
where $$\omega =\left( \begin{array}{c} \omega_{1}\\
\vdots \\
\omega_{n}
\end{array}\right)\in \R^n$$ is the frequencies vector. The 
associated system is the \textbf{harmonic oscillator} whose equations read
$$
\left\{ \begin{array}{cccc}
\dot q_{j}&=&\omega_{j} p_{j}, &\quad j=1,\ldots,n  \\
\dot p_{j}&=&-\omega_{j} q_{j}, &\quad j=1,\ldots,n 
\end{array}\right.
$$
and whose solutions are quasi-periodic functions given by
$$
\left\{ \begin{array}{cccc}
 q_{j}(t)&=&q_{j}(0) \cos \omega_{j}t + p_j({0}) \sin \omega_{j}t, &\quad j=1,\ldots,n  \\
 p_{j}(t)&=&-q_{j}(0) \sin \omega_{j}t + p_j({0}) \cos \omega_{j}t, 
 &\quad j=1,\ldots,n. 
\end{array}\right.
$$ 
Let us notice that for each $j$, $(q_{j},p_{j})$ describes a circle of 
radius\\ $\frac{p_{j}(0)^2+q_{j}(0)^2}{2}=:I_{j}$ and thus the orbits of the 
harmonic
oscillator are included in tori 
$$
 T_{I}:=\{(q,p)\in \R^{2n} \mid 
(p_{j}^2+q_{j}^2)/{2}=I_{j}, \ j=1,\ldots,n \}$$
whose dimension is generically $n$ (it 
can be less if ${p_{j}(0)^2+q_{j}(0)^2}=0$ for some $j$). To 
decide wether the orbit fills the torus or not we need the following 
definition:
\begin{definition}\label{nonresonant}
A frequencies vector, $\omega \in \R^n$, is \textbf{non resonant} if
 $$k\cdot \omega 
:=\sum_{j=1}^nk_{j}\omega_{j} \neq 0\quad \mbox{for all}\quad k\in \Z^n\setminus 
\{0\}.$$
\end{definition}
From number theory we learn that if $\omega$ is non resonant (or not 
rationnally dependent) then $\{ k\cdot \omega \mid k\in \Z^n \}$ is 
dense in $\R^n$ and thus we deduce that the orbit (or trajectory) is 
dense in the torus. On the contrary, if $\omega$ is resonant then the 
orbit is not dense in $T_{I}$ but in a torus of smaller dimension. 
For instance if all the frequencies are rationally proportional, 
$k_{1}\omega_{1}=k_{2}\omega_{2}=\ldots=k_{n}\omega_{n}$ for some 
choice of $k_{1},\ldots,k_{n}$ in $\Z$, the orbit is a circle and the 
solution is in fact periodic.

\subsection{Integrability}\label{int}
\begin{definition}
A \textbf{constant of motion}  (or an integral of motion)  for $H$ is a regular 
function, $F\in C^{\infty}(M,\R)$ satisfying $\{F,H\}=0$.
\end{definition}
\begin{proposition}
Let $F\in C^{\infty}(M,\R)$ then, if $t\mapsto (q(t),p(t))$ is a 
solution of the Hamiltonian system associated to $H$,
$$
\frac{d}{dt}F(q(t),p(t))=\{ F,H\} (q(t),p(t)).
$$
In particular, if $F$ is a constant of motion, then $F(q,p)$ is 
invariant under the flow generated by $H$.
\end{proposition}
\proof
By definition,
$$
\frac{d}{dt}F(q(t),p(t))= 
\sum_{j=1}^n \frac{\partial F}{\partial q_{j}}\dot q_{j}
+\frac{\partial F}{\partial p_{j}}\dot p_{j}= 
\sum_{j=1}^n \frac{\partial F}{\partial q_{j}}\frac{\partial H}{\partial 
p_{j}}
-\frac{\partial F}{\partial p_{j}}\frac{\partial H}{\partial q_{j}}
=\{ F,H\} .$$
\qed \\
In the case of the harmonic oscillator the actions $I_{j},\ 
j=1,\ldots,n$, defined by
$$
I_{j}=\frac{p_{j}^2+q_{j}^2}{2}
$$
are integrals of the motion : $\dot I_{j}=0$.
\begin{definition}\label{integrable}
A $2n$-dimensional Hamiltonian system is \textbf{integrable in the sense of  Liouville} if there exist $n$ 
regular 
functions $F_{1},F_{2},\ldots,F_{n}\ \in C^{\infty}(M,\R)$ such that
\begin{itemize}
\item[(i)] $\{F_{j},H\}=0$ for $j=1,\ldots,n$ (i.e. the $F_{j}$ are 
integrals of the motion).
\item[(ii)] $\{F_{j},F_{k}\}=0$ for $j,k=1,\ldots,n$ (i.e. the $F_{j}$ are 
in involution). 
\item[(iii)] $(\nabla_{q,p}F_{j})_{j=1,\ldots,n}$ are lineary 
independent.
\end{itemize}
\end{definition}
In view of this definition, the harmonic oscillator is  integrable in 
the sense of Liouville on  the dense open subset where none of the 
actions $I_{j}$ vanishes : it 
suffices to choose $F_{j}=I_{j}$.\\
However it is not Liouville integrable on the whole space: the 
actions are not independent everywhere. More generally,
a $2n$-dimensional Hamiltonian system which admits $n$ integrals in 
involution that are independent on a dense open subset is often called 
a \textbf{Birkhoff integrable Hamiltonian systems}.

\medskip

Examples of Liouville integrable system are obtained when the Hamiltonian 
depends only on $p$: $H(q,p)=h(p)$. In this case, often called the 
canonical example of integrable Hamiltonian system, $\dot p_{j}=0$ 
and $(p_{j})_{j=1,\ldots,n}$ satisfies (i), (ii) and (iii) of the 
previous definition. 
Actually the motion is trivial since $\dot q_{j}=\frac{\partial 
H}{\partial p_{j}}=$constant$=:\omega_{j}$ and thus we can \textit{integrate} 
the equations to obtain
$$
\left\{ \begin{array}{llll}
 q_{j}(t)&=&q_{j}(0)+ \omega_{j}t , &\quad j=1,\ldots,n  \\
 p_{j}(t)&=&p_{j}(0), &\quad j=1,\ldots,n. 
\end{array}\right.
$$

\medskip

Let $H$ be a Liouville integrable Hamiltonian  and denote by
$F_{1},F_{2},\ldots,F_{n}$ a complete set of independent integrals in 
involution on the phase space $M$. The leaves
$$
M_{c}=\{(q,p)\in M \mid F_{j}(q,p)=c_{j},\ j=1,\ldots,n\}
$$
are smooth submanifolds of $M$ of dimension and codimension $n$ 
\footnote{Actually they are Lagrangian submanifolds : 
submanifolds of 
maximal dimension such that the restriction of the symplectic form to 
it vanishes.}, and the whole manifold $M$ is foliated into these leaves.
\begin{theorem}
( Arnold-Liouville Theorem) Let
$H$ be a  Liouville integrable Hamiltonian on $M$ a symplectic 
manifold of dimension $2n$. If one of its leaves is 
compact and connected then there exists a neighborhood $\U$ of this 
leave, a 
neighborhood $D$ of $0$ in $\R^n$ and a change of 
variable\footnote{here $T^n=S^{1}\times\ldots\times S^1$, 
 $n$ times, is the $n$ dimensional torus}
 $\Psi : D\times T^n\ni(I,\theta)\mapsto 
(q,p)\in \U $
such that
\begin{itemize}
\item[(i)] $H\circ \Psi=h(I)$ is a function of the actions  alone
\item[(ii)] the Hamiltonian formalism is preserved, i.e., in the new 
variables, the equations read $\dot I_{j}=0,\ \dot \theta_{j}= -
\frac{\partial 
h}{\partial I_{j}},\ j=1,\ldots,n$ (i.e. the change of variable is a 
canonical transformation in the sense of the definition \ref{cano} 
below). 
\end{itemize}
\end{theorem}
That means that, every Liouville integrable Hamiltonian 
system with compact leaves is equivalent to the canonical one.
The new variables are called the 
\textbf{action-angle variables}.

In the case of the harmonic oscillator, the action-angle variables are 
given by the symplectic polar coordinates :
$$I_{j}=\frac{p_{j}^2+q_{j}^2}{2} ,\quad \theta_{j}=\arctan 
\frac{q_{j}}{p_{j}},$$
they are well defined on the dense open subset where none of the 
actions $I_{j}$ vanishes.

Notice that the Arnold-Liouville theorem implies that all the leaves 
$M_{c}$ are tori. Therefore, in this case,
the whole phase space $M$ is foliated by invariant tori of dimension 
$n$ (so called Lagrangian tori). This is not true in the case of a Birkhoff 
integrable Hamiltonian system where the dimension of the leaves may 
vary (as in the case of the harmonic oscillator).

\subsection{Perturbation of integrable Hamiltonian 
system}\label{perturb}
We consider a Hamiltonian function $H= H_{0}+P$ where $H_{0}$ is 
integrable and $P$ is a perturbation term. 

The general philosophy will consist in transforming $H$ in such a way 
that the new Hamiltonian 
be closer to an integrable one: $H\to \tilde H= \tilde H_{0} + \tilde 
P$ with $\tilde H_{0}$ still integrable and $\tilde P\ll P.$ The first 
question is: How to transform $H$? We cannot use all changes of 
variable because we want to conserve the Hamiltonian structure.
\begin{definition}\label{cano} A map $\varphi : M\ni (q,p)\mapsto (\xi, \eta)\in M$
    is a \textbf{canonical 
transformation} (or a \textbf{symplectic change of coordinates}) if 
\bi
\item $\varphi$ is a diffeomorphism
\item $\varphi$ preserves the Poisson Bracket : $\{F,G\}\circ \varphi=\{ 
F\circ \varphi,  G\circ\varphi \} for any *F$ and $G$.
\ei
\end{definition}
As a consequence, if $\tilde H= H\circ \varphi^{-1}$ with $\varphi$ 
canonical, then the Hamiltonian system reads in the new variables 
$(\xi,\eta)$ as in the old ones 
$$\dot \xi_{j}= \frac{\partial \tilde H}{\partial 
\eta_{j}}\ ,\quad \dot \eta_{j}= -\frac{\partial \tilde H}{\partial 
\xi_{j}}\quad j=1,\ldots,n.$$ 
There exists a very convenient way of constructing canonical 
transformation:\\
Let $\chi : M\to \R$ a regular function and denote $\varphi_{t}$ the 
flow generated by $X_{\chi}$. If $\varphi_{t}$ is well defined up to 
$t=1$, the map $\varphi :=\varphi_{1}$ is 
called the Lie transform associated to $\chi$. More explicitely, the 
new couple of variables $(\xi,\eta)= \varphi(q,p)$ is the value at 
time 1 of the solution of the system $\frac{d}{dt}\left(
\begin{array}{c}
    \xi  \\
    \eta  
\end{array}\right)=X_{\chi}(\xi,\eta)$ whose value at $t=0$ is $(q,p)$. 
Notice that, since the map $(t;q,p)\mapsto \varphi_{t}(q,p)$ is defined on an 
open set (cf. the Cauchy-Lipschitz theorem), if $\varphi_{1}$ is defined at 
the point $(q,p)$ then it 
is locally defined around $(q,p)$.
\begin{proposition}
A Lie transform is  canonical.
\end{proposition}
\proof : classical (see for instance \cite{A1}).

\medskip

The following lemma will be essential to use the Lie transforms 
\begin{lemma}\label{lie}
Let $G:M\to \R$ be a regular function. Then
$$\frac{d}{dt}(G\circ \varphi_{t})= \{G,\chi\}\circ \varphi_{t}.$$
\end{lemma}
\proof 
\begin{eqnarray*}
\frac{d}{dt}(G\circ \varphi_{t})(q,p)&=&\nabla G(\varphi_{t}(q,p))\cdot 
\dot\varphi_{t}(q,p)\\
&=& \nabla G(\varphi_{t}(q,p))\cdot 
J\nabla\chi(\varphi_{t}(q,p))\\
&=&\{G,\chi\}(\varphi_{t}(q,p)).
\end{eqnarray*}
\qed

Then, using the Taylor expansion of $G\circ \varphi_{t}$ at 
$t=0$, evaluated at $t=1$, we obtain for any $k\geq 0$,
$$G\circ \phi (q,p)=\sum_{l= 0}^k G_{l}(q,p) +(k+1)\int_{0}^1 
(1-t)^k G_{k+1}\circ \varphi_{t}(q,p) dt$$ with $G_{l}=1/l\{G_{l-1},\chi\}$ 
for $l\geq1$ and $G_{0}=G$.

\section{The Birkhoff normal form theorem in finite 
dimension}\label{section:birk-fini}
In this section we consider perturbations of Hamiltonian systems near 
an elliptic fixed point.

Let $H$ be a Hamiltonian function on $M$  having 
an isolated equilibrium. Without lost of generality we can assume that the 
equilibrium is at the origin in $\R^{2n}$, that the origin belongs 
to $M$ and that $H(0,0)=0$. Then the 
Hamiltonian reads
$$
H=\frac 1 2 <A(q,p),(q,p)>+ \mbox{ cubic terms}+\ldots
$$
where $A$ is the 
Hessian of $H$ at $0$, a symmetric $2n\times 2n$ real matrix. 
Since we suppose the equilibrium is elliptic, 
the spectrum of the linearized system $\dot{u}= JAu$ is purely 
imaginary:
$$\mbox{spec}(JA)= \{\pm i\omega_{1},\ldots,\pm i\omega_{n}\}$$
with $\omega_{1},\ldots\omega_{n}$ real. It turns out that there 
exists a linear symplectic change of coordinates that brings the 
quadratic part of $H$ into the following normal form (cf. \cite{HZ}, section1.7, 
theorem 8)
$$
<A(q,p),(q,p)>=\sum_{j=1}^n \omega_{j} ({p_{j}^2+q_{j}^2})
$$
where, for simplicity, we denote the new coordinates by the same 
symbols.

Therfore, in this section, we will focus on the perturbation of
the harmonic oscillator
\be\label{H0}
H_{0}(q,p)=\sum_{j=1}^n \omega_{j} \frac{p_{j}^2+q_{j}^2}{2}
=\sum_{j=1}^n \omega_{j} I_{j},\ee
where we denote  $I_{j}(q,p):=
\frac{p_{j}^2+q_{j}^2}{2}$ the $j^{th}$ action of $(q,p)$.\\
 The total
Hamiltonian reads $H=H_{0}+P$ where
$P$ is a regular real valued function, $P\in C^{\infty}(M,\R)$, which
is at least cubic, $P=O(\norma{(q,p)}^3)$, in 
such a way that, in a 
small neighborhood of $(0,0)$, $P$ will appear as a perturbation of 
$H_{0}$. 

We say that $P$ is in 
\textbf{normal form} with respect to $H_{0}$ if it commutes with the 
integrable part: $\{P,H_{0}\}=0$.\\
For $k\in \Z^n$, we denote by $\vak$ the length of $k$: $\vak = 
|k_{1}|+\ldots+|k_{n}|$. We need a 
refined version of the nonresonancy definition (compare with definition 
\ref{nonresonant}):
\begin{definition}\label{NR-faible}
Let $r\in \N$. A frequencies vector, $\omega \in \R^n$, is \textbf{non resonant up to 
order $r$} if
$$k\cdot \omega 
:=\sum_{j=1}^nk_{j}\omega_{j} \neq 0\quad \mbox{for all}\quad k\in 
\Z^n \mbox{ with } 0< \vak \leq r.$$
\end{definition}
Of course if $\omega$ is non resonant then it is nonresonant up to 
any order $r\in \N$.

\subsection{The theorem and its dynamical consequences}
We begin stating the classical Birkhoff normal form theorem 
(see for instance \cite{Mos68, HZ}).
\begin{theorem}\label{birk-fini} \textrm{\textbf{[Birkhoff Normal Form Theorem]}}
Let $H=H_{0}+P$, $H_{0}$ being the harmonic oscillator \eqref{H0} 
and $P$ being a $C^{\infty}$ real valued function having a zero of 
order 3 at the origin and fix $r\geq 3$ an integer. There exists $\tau : \U\ni 
(q',p')\mapsto (q,p)\in \V$ a real analytic canonical transformation 
from a neighborhood of the origin to a neighborhood of the origin
 which puts $H$ in normal form up 
to order $r$ i.e.
$$H\circ \tau= H_{0}+Z+R$$
with
\bi
\item[(i)] $Z$ is a polynomial or order $r$ and is in \textbf{normal form},
i.e.: 
$\{Z,H_{0}\}=0$.
\item[(ii)] $R\in C^\infty (M,\R)$ and $R(q',p')=O(\norma{(q',p')})^{r+1}$.
\item[(iii)] $\tau$ is close to the identity: $\tau 
(q',p')=(q',p')+O(\norma{(q',p')})^2$.
\ei
In particular if $\omega$ is non resonant up to order $r$ then $Z$ depends only on the 
new actions: 
$Z=Z(I'_{1},\ldots,I'_{n})$ with $I'_{j}=\frac{(p'_{j})^{2}+(q'_{j})^2}{2}$.    
\end{theorem}
Before proving this theorem, we analyse its dynamical consequences in 
the non resonant case.
\begin{corollary}\label{coro1}
 Assume $\omega$ is non resonant. For each $r\geq 3$ there exists 
 $\e_{0}>0$ and $C>0$ such that if $\norma{(q_{0},p_{0})}=\e <\e_{0}$ the 
 solution $(q(t),p(t))$ of the Hamiltonian system associated to $H$ 
 which takes value $(q_{0},p_{0})$ at $t=0$ 
 satisfies    
$$\norma{(q_{t},p_{t})}\leq 2\e \quad \mbox{for }\va{t}\leq 
\frac{C}{\e^{r-1}}.$$
Furthermore for each $j=1,\ldots,n$
$$|I_{j}(t)-I_{j}(0)|\leq \e^3 \quad \mbox{for }\va{t}\leq 
\frac{C}{\e^{r-2}},$$
where $I_{j}(t)=I_{j}(q_{t},p_{t})$.
\end{corollary}
\proof 

Denote $z:=(q,p)$ and $z'=\tau^{-1}(z)$ where $\tau$ is the 
transformation given by theorem \ref{birk-fini}. Then we define
$N(z):=\norma{z}^2=2\sum_{j=1}^n I_{j}(q,p)$. Using that $Z$ depends 
only on the new actions, we have
$$
\{N,H\}(z)=\{N\circ \tau ,H\circ \tau \}\circ \tau^{-1}(z)= \{N\circ 
\tau, R\}(z')= O(\norma{z'}^{r+1})= O(\norma{z}^{r+1})$$ 
Therefore $|\dot N|\leq C N^{(r+1)/2 }$. Using that the solution to the 
ODE $\dot x=ax^d$ is given by ($d>1$)
$$
x(t) =x_{0}(1-x_{0}^{d-1}(d-1)at)^\frac{-1}{d-1}
$$
one easily deduces the first part of the corollary.

To prove the second part, write (with $I'=I\circ \tau^{-1}$)
$$
\va{I_{j}(t)-I_{j}(0)}\leq 
\va{I_{j}(t)-I'_{j}(t)}+\va{I'_{j}(t)-I'_{j}(0)}+\va{I'_{j}(0)-I_{j}(0)}.$$
The first and the third term of the right side of this inequality are 
estimated by $c\e^3$ because $\norma{z-z'}\leq c\norma{z}^2$ and 
$|I'-I|\leq \norma{z'-z}\norma{z+z'}$. To estimate the middle term we 
notice that
\be \label{estimI}
\frac{d}{dt}I'_{j}=\{ I'_{j},H \}=\{ I_{j},H\circ \tau \}\circ \tau^{-1} = O(\norma{z'}^{r+1})
\ee
and therefore, for $\vat\leq \frac{c}{\e^{r-2}}$,
$$
\va{I'_{j}(t)-I'_{j}(0)}\leq c'\e^3.
$$\qed

We can also prove that the solution remains close to a torus for a 
long time, that is the contain of the following
\begin{corollary}\label{coro2}
 Assume $\omega$ is non resonant. For each $r\geq 3$ there exists 
 $\e_{0}>0$ and $C>0$ such that if $\norma{(q_{0},p_{0})}=\e <\e_{0}$ 
 then there exists a torus $\mathcal T_{0}\subset M$ satisfying
$$\mbox{dist}((q(t),p(t)),\mathcal T_{0})\leq C\e^{r_{1}} \quad 
\mbox{for } \vat \leq 1/\e^{r_{2}}$$
where $r_{1}+r_{2}=r+1$.
\end{corollary}
\proof 
Let $$\mathcal T_{0}:=\{ (q,p)\mid I_{j}(\tau^{-1}(q,p))=
I_{j}(\tau^{-1}(q_{0},p_{0})),\ j=1,\ldots,n\}.$$
Using \eqref{estimI} we deduce that for $\vat \leq 1/\e^{r_{2}}$
$$
\va{I'_{j}(t)-I'_{j}(0)}\leq c'\e^{r_{1}},
$$
where as before $I'=I\circ \tau^{-1}$.
Therefore, using assertion (iii) of theorem \ref{birk-fini} we obtain
the thesis. \qed

\begin{remark} \textbf{An extension to the Nekhoroshev's theorem} \\
If $\omega$ is non resonant at any order, it is natural to try to optimize the 
choice of $r$ in theorem \ref{birk-fini} or its corollaries. 
Actually, if we assume that $\omega$ satisfies a diophantine condition
$$
|k\cdot \omega|\geq \gamma \vak^{-\alpha} \mbox{ for all } k\in 
\Z^n\setminus \{0\},$$
standart estimates (see for instance \cite{ben85, GG85, P93}) allow to 
prove that for $(q,p)$ in $B_{\rho}$, the ball centered at the origin and of 
radius $\rho$, the remainder term in theorem \ref{birk-fini} is of 
order $(r!)^{\tau +1}\rho^{r+1}$. This leads to show that the 
constant $C$ in corollary \ref{coro1} is of order $(r!)^{-(\tau +1)}$. 
Namely one proves that if $\norma{(q_{0},p_{0})}=\e$ is small enough
$$\norma{(q_{t},p_{t})}\leq 2\e \quad \mbox{for }\va{t}\leq 
\frac{C }{\e^{r-1}(r!)^{\tau +1}}$$
where the new constant $C$  depends only on $P$ and on the number of 
degrees of 
freedom\footnote{This dependence with respect to $n$ makes impossible 
to generalize, at least easily, this remark in the infinite 
dimensional case.} $n$. Using the Stirling's formula for $r!$ and 
choosing $r=e^2 \e^{-1/(\tau +1)}$, one obtains that the solution 
remains controled by $2\e$ during an \textbf{exponentially long time}: 
$$\norma{(q_{t},p_{t})}\leq 2\e \quad \mbox{for }\va{t}\leq C 
\exp \left(\frac{\beta}{\e^{1/(\tau +1)}}\right),$$
where $\beta$ is a non negative constant. This last statement is a 
formulation of the Nekhoroshev's theorem (see \cite{Nek77}).
\end{remark}    
\subsection{Proof of the Birkhoff normal form theorem}

We prefer to use the complex variables
$$\xi_{j}=\frac{1}{\sqrt 2}(q_{j}+ip_{j}), \quad 
\eta_{j}=\frac{1}{\sqrt 2}(q_{j}-ip_{j})$$
because the calculus are easier in this framework. Notice in 
particular that the actions read $I_{j}=\xi_{j}\eta_{j}$ and thus it 
is very simple to express that a monomial 
$\xi_{j_{1}}\ldots\xi_{j_{k}}\eta_{l_{1}}\ldots\eta_{l_{k'}}$ depends 
only on the actions, it suffices that $k=k'$ and 
$\{ j_{1},\ldots,j_{k}\}=\{l_{1},\ldots,l_{k}\}$.

We have $H_{0}=\sum_{j=1}^n \omega_{j}\xi_{j}\eta_{j}$ and
we easily verify that, in these variables, the Poisson bracket 
reads
$$
\{ F,G \}=i\sum_{j=1}^n \ \frac{\partial F}{\partial \xi_{j}}
\frac{\partial G}{\partial \eta_{j}}-
\frac{\partial F}{\partial \eta_{j}}\frac{\partial G}{\partial \xi_{j}}
\ .$$
We will say that a function $F$ defined in the variable $(\xi, \eta)$ 
 is \textbf{real} when $F(\xi,\bar \xi)$ is real which means that in 
 the original variables $(q,p)$, $F$ is real valued.

We now begin the proof of theorem \ref{birk-fini}. Having fixed some
$r\geq 3$, the idea is to construct iteratively for $ k=
2,\ldots,r$,  a canonical transformation $\tau_k$, defined on a 
neighborhood of the origin in $M$, 
and real functions $Z_{k}, P_{k+1}, R_{k+2}$ such that
\begin{equation} 
\label{b.1} 
H_{k}:=H\circ\tau_k=H_{0}+Z_{k}+P_{k+1}+R_{k+2} 
\end{equation} 
and with the following properties
\bi
\item[(i)] $Z_{k}$ is a polynomial of degree $k$ having a zero of 
degree 3 at the origin and $\{Z_{k},H_{0}\}=0$.
\item[(ii)] $P_{k+1}$ is a homogeneous polynomial of degree $k+1$.
\item[(iii)] $R_{k+2}$ is a regular Hamiltonian having a zero of 
order $k+2$ at the origin.
\ei
Then \eqref{b.1} at order $r$ proves theorem \ref{birk-fini} with 
$Z=Z_{r}$ and $R=P_{r+1}+R_{r+2}$.

First remark that the Hamiltonian $H=H_{0}+P$ has the form 
(\ref{b.1}) with $k=2$, $\tau_2=I$, $Z_{2}=0$, 
$P_{3}$ being the Taylor's polynomial of $P$ at degree 
$3$ and $R_{4}=P-P_{3}$. We 
show now how to pass from $k$ to $k+1$.\\
We search for $\tau_{k+1}$ of the form $\tau_{k}\circ \phi_{k+1}$, 
$\phi_{k+1}$ being a Lie transform associated to the Hamiltonian 
function $\chi_{k+1}$. Recall from section \ref{perturb} that for 
regular $F$
$$
F\circ \phi_{k+1}= F
+\{F,\chi_{k+1}\}+1/2\{\{F,\chi_{k+1}\},\chi_{k+1}\} +\ldots$$
We search for $\chi_{k+1}$ as a homogeneous real polynomial of 
degree $k+1$ and we decompose $H_{k}\circ \phi_{k+1}$ as follows
\begin{eqnarray} 
\label{b.6} 
H_{k}\circ \phi_{k+1}&=& H_{0}+Z_{k}+\{H_{0},\chi_{k+1}\}+P_{k+1} 
\\ 
\label{b.7} 
&+& R_{k+2}\circ \phi_{k+1}\ +\ H_{0}\circ \phi_{k+1}-H_{0}-\{H_{0},\chi_{k+1}\}
\\ 
\label{b.8} 
&+& Z_{k}\circ \phi_{k+1} -Z_{k}\
+\ P_{k+1}\circ \phi_{k+1}-P_{k+1}\ .
\end{eqnarray} 
Notice that if $F_{1}$ is a homogeneous polynomial of degree 
$d_{1}$ and $F_{2}$ is a homogeneous polynomial of degree 
$d_{2}$ then $\{F_{1},F_{2}\}$ is a homogeneous polynomial of degree 
$d_{1}+d_{2}-2$. Notice also that, since 
$\chi_{k+1}(\xi,\eta)=O(\norma{(\xi,\eta)})^{k+1}$, we have 
$$\phi_{k+1}(\xi,\eta)=(\xi,\eta)+O(\norma{(\xi,\eta)})^{k}.$$
Using 
these two facts we deduce that \eqref{b.7} and \eqref{b.8} are  regular 
Hamiltonians 
having a zero of order $k+2$ at the origin. Therefore, using the 
Taylor formula, the sum of these terms decomposes in $P_{k+2}+R_{k+3}$ with $P_{k+2}$ 
and $R_{k+3}$ satisfying the properties (ii) and (iii). So it remains 
to prove that $\chi_{k+1}$ can be choosen in such a way that
$Z_{k+1}:=Z_{k}+\{H_{0},\chi_{k+1}\}+P_{k+1}$ 
satisfies (i). This is a consequence of the following lemma
\begin{lemma}\label{homolo}
Let $Q$ be a homogeneous real  polynomial of degree $k$, there exist 
two
homogeneous real valued polynomials $\chi$ and $Z$ of degree $k$ such that
\be \label{homo}\{H_{0},\chi\}+Q=Z\ee
and
\be \label{normal}\{Z,H_{0}\}=0.\ee

\end{lemma}
Equation \eqref{homo} is known in the literature as the 
\textbf{homological equation}.
\proof
For $j\in [1,\ldots,n]^{k_{1}}$ and $l\in [1,\ldots,n]^{k_{2}}$, 
denote $\xi^{(j)}= \xi_{j_{1}}\ldots
\xi_{j_{k_{1}}}$ and $\eta^{(l)}=\eta_{l_{1}}\ldots\eta_{l_{k_{2}}}$.
A direct calculus shows that 
$$
\{H_{0}, \xi^{(j)}\eta^{(l)}
\}
=-i\Omega(j,l) \xi^{(j)}\eta^{(l)}$$
with
$$
\Omega(j,l):= 
\omega_{j_{1}}+\ldots+\omega_{j_{k_{1}}}-\omega_{l_{1}}-\ldots 
-\omega_{l_{k_{2}}}.$$
Let
$$
Q=
\sum_{(j,l)\in [1,\ldots,n]^{k}}a_{jl}\xi^{(j)}\eta^{(l)}
$$
where $(j,l)\in [1,\ldots,n]^{k}$ means that $j\in 
[1,\ldots,n]^{k_{1}}$ and $l\in [1,\ldots,n]^{k_{2}}$ with $k_{1}+k_{2}=k$.
Then defining
$$b_{jl}=i\Omega(j,l)^{-1}a_{ij},\quad c_{jl}=0\quad \mbox{when }
\Omega(j,l)\neq0$$
and
$$c_{jl}=a_{ij},\quad b_{jl}=0\quad \mbox{when }
\Omega(j,l)=0\, ,$$
the polynomials
$$
\chi=\sum_{(j,l)\in [1,\ldots,n]^{k}}b_{j,l}\xi^{(j)}\eta^{(l)}
$$
and
$$
Z=\sum_{(j,l)\in [1,\ldots,n]^{k}}c_{j,l}\xi^{(j)}\eta^{(l)}
$$
satisfy \eqref{homo} and \eqref{normal}. Furthermore, that $Q$ is real 
is a consequence of the symmetry relation: $\bar a_{jl}=a_{lj}$. Taking 
into acount that $\Omega_{lj}=-\Omega_{jl}$,
this symmetry remains satisfied for the polynomials 
$\chi$ and $Z$. \qed

To complete the proof of theorem \ref{birk-fini}, it remains to 
consider the non resonant case. Recall that we use lemma \ref{homolo}
to remove succesively  parts of the polynomials $P_{k}$ for 
$k=3,\ldots,r$. Therefore the $\Omega_{j,l}$ that we need to consider can be 
written $k\cdot \omega$ for a $k\in \Z^n$ satisfying $\vak \leq r$. 
Thus if $\omega$ is nonresonant up to order $r$,
these $\Omega_{j,l}$ can 
vanish only if $j=l$ and thus the normal terms constructed in lemma 
\ref{homolo} have the form $Z=\sum_{j}a_{j,j}\xi^{(j)}\eta^{(j)}=
\sum_{j}a_{j,j}I^{(j)}$, i.e. $Z$ depends only on the actions. \qed

\begin{exercise} Let $Q=\xi_{1}\eta_{2}^2$ (resp. $Q=\xi_{1}\eta_{2}^2+
    \xi_{1}^2\eta_{2}$) and assume $\omega_{1}/\omega_{2}\notin \mathbb 
    Q$. Compute the corresponding $\chi$, $Z$. Then compute the new 
    variables $(\xi',\eta')=\tau^{-1}(\xi,\eta)$, $\tau$ being the 
    Lie transform generated by $\chi$. Verify that 
    $H_{0}(\xi',\eta')= H_{0}(\xi,\eta) +Q(\xi,\eta)$ (resp.
    $H_{0}(\xi',\eta')= H_{0}(\xi,\eta) +Q(\xi,\eta) +\mbox{ order }4$).
\end{exercise}

\section{A Birkhoff normal form theorem in infinite dimension }\label{sect-infini}
In this section we want to generalize the Birkhoff normal form theorem stated and 
proved in section \ref{section:birk-fini} in finite dimension to the 
case of infinite dimension. In view of section 
\ref{section:birk-fini}, we can previse that we face to difficulty: 
first we have to replace the definition \ref{nonresonant} by a concept 
that makes sense in infinite dimension. This 
will be done in definition \ref{SNR}. The second difficulty is more 
structural: we have to define a class of perturbations $P$, and in 
particular a class of polynomials, in which the Birkhoff procedure can 
apply even with an infinite number of monomials. Concretly, the problem is to be 
able to verify that at each step of the procedure the formal polynomials 
that we construct are at least continuous function on the phase 
space (the continuity of a polynomial is not automatic in infinite 
dimension since they may contain an infinite number of monomials). 
This class of polynomial is defined in definition \ref{poly} and is directly 
inspired by a class of multilinear forms introduced in \cite{DS1}, 
\cite{DS2}. In section \ref{model} we define our model of infinite 
dimensional integrable Hamiltonian system  and in section \ref{thm} we state the 
Birkhoff type result and its dynamical consequences.

\subsection{The model}\label{model}
To begin with we give an abstract model of infinite dimensional 
Hamiltonian system. In section \ref{ex} we will give some concrete  PDEs 
that can be described in this abstract framework.

We work in the phase space $\Ps\equiv\Ps(\R):=l^2_{s}(\R)\times l^2_{s}(\R)$
where, for $s\in \R$,  
$l^2_{s}(\R):=\{(a_{j})_{j\geq 1}\in \R^\N \mid \sum_{j\geq 
1}j^{2s}\va{a_{j}}^2   \}$  is a Hilbert space for the 
standart norm: $\norma{a}_{s}^2=\sum_{j\geq 
1}\vaj^{2s}\va{a_{j}}^2$. 

Let us denote by $(\cdot ,\cdot )$ the 
$l^2$-scalar product on $l^2_{s}(\R)$.
Let $\U$ be an open subset of $l^2_{s}(\R)$, for $F\in C^{1}(\U,\R)$
and $a\in \U$, we define the $l^2$ gradient $\nabla F(a)$ by
$$
DF(a)\cdot h =(\nabla F(a),h), \quad \mbox{for all }h\in 
l^2_{s}(\R)$$
where $DF(a)$ denotes the differential of 
$F$ at the point $a$. We write $$\nabla F(a)\equiv \left( 
\frac{\partial F}{\partial a_{j}}(a)\right)_{j\geq 1}.$$ 
Notice that, without further hypothesis on $F$, we only 
have $\nabla F(a)\in l^2_{-s}(\R)$.\\
Then we endow $\Ps$ with the canonical 
symplectism $\sum_{j\geq 1}dq_{j}\wedge dp_{j} $ and we 
define the Hamiltonian vector field of a regular Hamiltonian function 
on an open subset $\U$ of $\Ps$,  
$H \in  C^{\infty}(\U,\R)$ by
$$
X_{H}(q,p)=\left(\begin{array}{c}\left( 
\frac{\partial H}{\partial p_{j}}(q,p)\right)_{j\geq 1}\\
-\left( 
\frac{\partial H}{\partial q_{j}}(q,p)\right)_{j\geq 1}\end{array}\right).
$$
Again without further hypothesis on $H$, we only 
know that  $X_{H}(q,p)\in l^2_{-s}(\R)\times l^2_{-s}(\R)$. However, in 
order to consider the flow of the Hamilton's equations
$$
\frac{d}{dt}
\left(\begin{array}{c}
    q  \\
    p
\end{array}\right)=X_{H}(q,p)$$
on $\Ps$, we prefer that the vector field preserves this phase space\footnote{
this condition is not 
really necessary, $X_{H}(q,p)$ could be unbounded as an operator from 
$\Ps$ to $\Ps$ }, i.e. 
$X_{H}(q,p)\in \Ps$ for $(q,p)\in \Ps$. Thus we will be interested in the following class 
of admissible Hamiltonian functions
\begin{definition}
Let $s\geq 0$, we denote by 
$\Hs$ the space of real valued  functions $H$ defined on
a neighborhood $\U$ of the origin in $\Ps$ and
satisfying 
$$H\in C^\infty (\U, \R)\quad \mbox{and} \quad X_{H}\in C^\infty (\U, \Ps).$$
\end{definition}
In particular the Hamiltonian vector fields of functions $F,\ G$ in 
$\Hs$ are in $l^2_{s}(\R)\times l^2_{s}(\R)$ and we can define their Poisson 
bracket by
$$\{ F,G \}(q,p)=\sum_{j\geq 1} \ \frac{\partial F}{\partial q_{j}}(q,p)
\frac{\partial G}{\partial p_{j}}(q,p) -
\frac{\partial F}{\partial p_{j}}(q,p)\frac{\partial G}{\partial q_{j}}(q,p)
\ .$$
We will also use the complex variables 
$$\xi_{j}=\frac{1}{\sqrt 2}(q_{j}+ip_{j}), \quad 
\eta_{j}=\frac{1}{\sqrt 2}(q_{j}-ip_{j}).$$ We have $(\xi, 
\eta)\in \Ps(\C)$, the complexification of $\Ps(\R)$.
In these variables, the Poisson bracket of two functions in 
$\Hs$
reads
$$
\{ F,G \}=i\sum_{j\geq 1} \ \frac{\partial F}{\partial \xi_{j}}
\frac{\partial G}{\partial \eta_{j}}-
\frac{\partial F}{\partial \eta_{j}}\frac{\partial G}{\partial \xi_{j}}
\ $$
where
$$\frac{\partial}{\partial \xi_{j}}=\frac{1}{\sqrt 2}\left(
\frac{\partial}{\partial q_{j}}-i\frac{\partial}{\partial 
p_{j}}\right), \quad 
\frac{\partial}{\partial \eta_{j}}=\frac{1}{\sqrt 2}\left(
\frac{\partial}{\partial q_{j}}+i\frac{\partial}{\partial 
p_{j}}\right).$$

As in the finite dimensional case, we will say that a function $F$ 
defined in the variable $(\xi, \eta)$ 
 is real when $F(\xi,\bar \xi)$ is real which means that in 
 the original real variables $(q,p)$, $F$ is real valued. We will use the 
 notation $z=(\ldots,\xi_{2},\xi_{1},\eta_{1},\eta_{2},\ldots)\in 
 l^2_{s}(\Zb, \C)$ where $\Zb=\Z\setminus\{0\}$. We will also 
 denote $\Nb =\N \setminus\{0\}$.

\medskip

Our model of integrable system is the harmonic oscillator
$$H_{0}=\sum_{j\geq 1} \omega_{j}\xi_{j}\eta_{j}$$
where $\omega =(\omega_{j})_{j\geq 1}\in \R^\N$ is the frequencies 
vector. We will assume that these frequencies grow at most polynomilally, 
i.e. that there exist $C>0$ and  $d\geq 0$ such that for any $j\in \Nb$,
 \be \label{om}|\omega_{j}|\leq C \vaj^d,\ee
 in such a 
way that  $H_{0}$ be well defined on $\Ps$ for $s$ large enough.
The perturbation term is a real  function, $P\in 
\Hs$, having a zero of order at least $3$ at the 
origin. Our Hamiltonian function is then given by
$$H=H_{0}+P$$
and the Hamilton's equations read, in the real variables,
\be \label{HeqR}
\left\{\begin{array}{ccc}
\dot q_{j}=&\omega_{j} p_{j} +\frac{\partial P}{\partial p_{j}},& 
j\geq 1\\
\dot p_{j} =& -\omega_{j}q_{j}-\frac{\partial P}{\partial q_{j}},& 
j\geq 1
\end{array}
\right. \ee
and in the complex ones
\be \label{HeqC}
\left\{\begin{array}{ccc}
\dot \xi_{j} =-&i\omega_{j} \xi_{j} -i\frac{\partial P}{\partial \eta_{j}},& 
j\geq 1\\
\dot \eta_{j} =& i\omega_{j}\eta_{j}+i\frac{\partial P}{\partial \xi_{j}},& 
j\geq 1.
\end{array}
\right. 
\ee

Our theorem will require essentially two hypotheses: one on the 
perturbation $P$ (see definition \ref{po} below) and one on the frequencies vector 
$\omega$ (see definition \ref{SNR} below). We begin by giving a motivation for 
these intriguing definitions.

As in the finite dimensional case, the game will consist in removing iteratively, 
by a canonical transform,
the cubic terms of $P$ that are not in normal form with respect to 
$H_{0}$, then the quartic ones and so on. The basic rule remains the 
same:
to remove the monomial $a_{j}z_{j_{1}}\ldots z_{j_{k}}$, we have to 
control the monomial $\frac{a_{j}}{\Omega(j)}z_{j_{1}}\ldots 
z_{j_{k}}$ where, as in the finite dimensional case (cf. proof of 
lemma \ref{homolo}), $\Omega(j)= \mbox{sign}(j_{1})\omega_{j_{1}}+\ldots+
\mbox{sign}(j_{k})\omega_{j_{}}$ is the small divisor. In contrast with the finite 
dimensional case, the number of monomials that we have to remove at 
each step is, a priori, infinite. Fortunately, the vector field of many of them 
are already small in the $l^2_{s}$-norm for $s$ large enough: \\
Consider the simple case where $P$ is a monomial of degree $k$, $P= z_{j_{1}}\ldots z_{j_{k}}$. 
Assume that $|j_{1}|\leq\ldots\leq |j_{k}|$ and that the three largest 
indexes are large and of the same order, say $N\leq |j_{k-2}|\leq 
|j_{k-1}|\leq |j_{k}|\leq 2N$. Then one gets
\begin{eqnarray*}
\norma{X_{P}(z)}_{s}^2 &=& \sum_{l\in \Zb} l^{2s}
\Va{\frac{\partial P}{\partial z_{l}}(z)}^2 \\
&\leq & \sum_{l=1}^k |j_{l}|^{2s}\va{z_{j_{1}}}^2\ldots\va{z_{j_{l-1}}}^2
\va{z_{j_{l+1}}}^2\ldots\va{z_{j_{k}}}^2 \\
&\leq & \sum_{l=1}^k \frac{|j_{l}|^{2s}}{\prod_{m\neq l}|j_{m}|^{2s}}
\prod_{m\neq l}|j_{m}|^{2s}\va{z_{j_{m}}}^2 \\
&\leq &  \frac{C}{N^{2s}}\norma{z}_{s}^{2k-2}
\end{eqnarray*}
which is small for large $N$.
This calculus explains why we will   
control only the small divisors that involve at most two large indexes, 
whence the definition \ref{SNR} below. Concerning the class of 
nonliniearities $P$ that we can consider, the preceding calculus does 
not suffice to justify the precise form of definition \ref{po} below 
but again it explains why the estimate \eqref{a} involves only the three 
largest indexes. Actually some 
other constraints are required like the control of the 
$l^{2}_{s}$-norm of $X_{P}$ (cf. proposition \ref{tame}) and the 
stability of the class under Poisson bracketing (cf. proposition 
\ref{bra}).

\medskip

For $j\in \Zb^k$ with $k\geq 3$, we define
$\mu(j)$ as the third largest integer between 
$\va{j_{1}},\ldots,\va{j_{k}}$. Then we set 
$S(j):={\va{j_{i_{0}}}-\va{j_{i_{1}}}} +\mu(j)$ where 
$\va{j_{i_{0}}}$ and $\va{j_{i_{1}}}$ are respectively the largest integer and the second 
largest integer between 
$\va{j_{1}},\ldots,\va{j_{k}}$. In particular
$$ \mbox{if } 
\va{j_{1}}\leq\ldots\leq\va{j_{k}} \mbox{ then }
\mu(j):=\va{j_{k-2}} \mbox{ and } 
S(j)=\va{j_{k}}-\va{j_{k-1}}+\va{j_{k-2}}.$$
For $j\in \Zb^k$ with $k\leq 2$, we fix $\mu(j)=S(j)=1$. 
\begin{definition}\label{po} Let $k\geq 3$, $N\in \N$ and $\nu \in [0,+\infty)$ and
let $$Q(\xi, \eta) \equiv Q(z)=\sum_{l=0}^k\sum_{j\in 
\Zb^l}a_{j}z_{j_{1}}\ldots z_{j_{l}}$$ be a formal polynomial of degree $k$
on $\Ps(\C)$. $Q$ is in the class 
$\T_{k}^{N,\nu}$ if there exists a constant $C>0$ such that for all 
$j$
\be\label{a}\va{a_{j}}\leq C \frac{\mu(j)^{N+\nu}}{S(j)^N}.\ee
\end{definition}
We will see in section \ref{poly} that $\T_{k}^{N,\nu}\subset \Hs$ 
for $s\geq \nu +1/2$ (cf. proposition \ref{tame}) and thus in 
particular a polynomial in $\T_{k}^{N,\nu}$ is well defined on a  
neighborhood of the origin in $\Ps(\C)$ for $s$ large enough.
The best constant $C$ in \eqref{a} defines a norm for which 
$\T_{k}^{N,\nu}$ is a Banach space. We set $$\T_{k}^{\infty,\nu}=\cap_{N\in 
\N}\T_{k}^{N,\nu}$$
and
$$\T^{\nu}= \cup_{k\geq 0}\T_{k}^{\infty,\nu}.$$
This definition is similar to a class of multilinear 
forms first introduced by Delort and Szeftel in \cite{DS1} and 
\cite{DS2}.

\begin{definition}\label{T} A function $P$ is in the class 
$\T$ if
\begin{itemize}
    \item
there exist $s_{0}\geq 0$ such that, for any $s\geq s_{0}$,  
$P\in \Hs$  
\item
for each $k\geq 1$ there exits $\nu \geq 0$ such that the Taylor's 
expansion  of 
degree $k$ of $P$ at zero belongs 
to $\T_{k}^{\infty,\nu}$.
\end{itemize}
\end{definition}
In  section \ref{poly} we will establish some properties of polynomials in 
$\T_{k}^{N,\nu}$, in particular we will see that this class has 
a good behaviour regarding to 
 the Poisson bracket (cf. proposition 
\ref{bra}).

\medskip

Concerning the frequencies, we define:
\begin{definition}\label{SNR} A frequencies vector $\omega \in \R^\Nb$ 
is \textbf{strongly non resonant} if for
    any $r\in\Nb$, there are $ \gamma >0$ and 
    $\alpha >0$ such that for any $j\in \Nb^{r}$ and any $1\leq i\leq 
    r$, one has 
    \begin{equation} 
    \label{A.2} 
    \left|\omega_{j_1}+\cdots+\omega_{j_{i}}-\omega_{j_{i+1}}-\cdots 
    -\omega_{j_{r}} \right|\geq \frac{\gamma}{\mu (j)^{\alpha}}  
    \end{equation} 
except if $\{j_{1},\ldots,j_{i}\}=\{j_{i+1},\ldots,j_{r}\}$.    
\end{definition}

This definition was first introduced  
in \cite{Bam03}.

\begin{remark}
The direct generalization of the definition \ref{nonresonant} to the 
infinite dimensional case would read
$\sum_{j\geq 1}\omega_{j}k_{j}\neq 0$ for all $k\in \Z^\Nb \setminus 0$. 
But in the infinite dimensional case this condition no more implies 
that there exists $C(r)$ such that \be \label{non}
\Va{\sum_{j\geq 1}\omega_{j}k_{j}}\geq C(r) 
\mbox{ for all } k\in \Z^\Nb \setminus 
0 \mbox{ satisfying } \sum_{j\geq 1}\va{k_{j}}\leq r\ee
which is the property 
that we used in the proof of theorem \ref{birk-fini}. Actually this 
last property \eqref{non} is too restrictive in infinite dimension. For 
instance when the frequencies are the eigenvalues of the 1-d 
Schr\"odinger operator with Dirichlet boundary conditions (cf. example
 \ref{ex.5}), one shows that  $\omega_{j}= j^2+ a_{j}$ where 
$(a_{j})_{j\geq 1}\in l^2$ (cf. for instance \cite{Mar86, PT}) and thus 
if $l$ is an odd integer and 
$j=(l^2-1)/2$ then one has 
$\omega_{j+1}-\omega_{j}-\omega_{l}\to_{l\to \infty}0$.
\end{remark} 

Our strongly nonresonant condition says that 
$\va{\sum_{j\geq 1}\omega_{j}k_{j}}$ is controled from below by a quantity which 
goes to zero when the third largest index of the frequencies involved grows to 
infinity, but the length of $k$ is fixed.  Precisely one has:
\begin{proposition}\label{SNRbis} A frequencies vector $\omega \in \R^\Nb$ 
is {strongly non resonant} if and only if for
any $r\in\Nb$, there are $ \gamma >0$ and $\alpha >0$ such that for 
any $N\in \Nb$
\be \label{A.2bis}\va{\sum_{m=1}^N\omega_{m}k_{m} 
+k_{l_{1}}\omega_{l_{1}}+k_{l_{2}}\omega_{l_{2}}}\geq \frac \gamma {N^\alpha}
\ee
for any indexes $l_{1},l_{2}>N$ and for any $k\in \Zb^{N+2}\setminus 
\{0\}$ with $\sum_{m=1}^N |k_{m}| \leq r$,  
 $|k_{l_{1}}|+|k_{l_{2}}|\leq 2$.
\end{proposition}
   In this form, the strongly nonresonant condition may be compared 
   to the so called Melnikov 
condition used in the KAM theory (cf.  
section \ref{KAM}).

\proof 
In order to see that the first form implies the second ones, we 
remark that the expression $\sum_{m=1}^N\omega_{m}k_{m} 
+k_{l_{1}}\omega_{l_{1}}+k_{l_{2}}\omega_{l_{2}}$ may be rewrite as
$\omega_{j_1}+\cdots+\omega_{j_{i}}-\omega_{j_{i+1}}-\cdots 
    -\omega_{j_{r'}}$ for some $r\leq r'\leq r+2$ and for some $j\in 
    \Nb^{r'}$ satisfying $\mu(j)\leq N$.\\
Conversely,  $\omega_{j_1}+\cdots+\omega_{j_{i}}-\omega_{j_{i+1}}-\cdots 
    -\omega_{j_{r}}$ may be rewrite as  
    $\sum_{m=1}^N\omega_{m}k_{m} 
    +k_{l_{1}}\omega_{l_{1}}+k_{l_{2}}\omega_{l_{2}}$ with 
    $N=\mu(j)$, $k_{l_{1}}, k_{l_{2}}=\pm 1$ and $\sum_{m=1}^N |k_{m}| 
    \leq r-2$.
\qed

\subsection{The statement}\label{thm}
We can now state our principal result: 
    
\begin{theorem}\label{birk-infini}
    Assume that $P$ belongs to  the class $\T$   and 
    that $\omega$ is strongly non resonant and satisfies \eqref{om} for 
    some $d\geq 0$. Then 
    for any $r\geq 3$ there exists $s_{0}$ and  for any $s\geq 
    s_{0}$ there exists $\U_{s}$, $\V_{s}$ neighborhoods of the origin in 
    $\Ps(\R)$ and 
    $\tau_{s} : \V_{s} \to \U_{s}$ a real analytic canonical transformation 
    which is the restriction to $\V_{s}$ of $\tau := \tau_{s_{0}}$ and  
     which puts $H=H_{0}+P$ in normal form up 
    to order $r$ i.e.
    $$H\circ \tau= H_{0}+Z+R$$
    with
    \bi
    \item[(i)] $Z$ is a continuous polynomial of degree $r$ with a 
    regular vector field (i.e. $Z\in \mathcal H^t$ for all $t\geq 0$) which only 
    depends on the actions: $Z= Z(I)$.
    \item[(ii)] $R\in \Hs(\V_{s},\R)$ and $\norma{X_{R}(q,p)}_{s}\leq 
    C_{s}\norma{(q,p)}_{s}^{r}$ for all $(q,p)\in \V_{s}$.
    \item[(iii)] $\tau$ is close to the identity: $\norma{\tau 
    (q,p)-(q,p)}_{s}\leq C_{s}\norma{(q,p)}_{s}^2$ for all $(q,p)\in \V_{s}$.
    \ei
\end{theorem}
This theorem was first proved in \cite{BG04} under a slightly more general 
hypothesis on the perturbation (cf. remark \ref{rktame} and 
\ref{rem:tame}).

This theorem says, as in the finite dimensional case, that we can 
change the coordinates in a neighborhood of the origin
in such a way that the Hamiltonian be 
integrable up to order $r$, $r$ being fixed at the principle. Remark 
that the concept of integrability that we gave in definition \ref{integrable} 
does not directly extend to the infinite dimensional case\footnote{It 
can be done with an appropriate definition of linear independence of 
an infinity of vector fields.}. However, if a 
Hamiltonian $H(q,p)$ depends only on the actions $I_{j},$ $j\geq 1$ 
then we can say that $H$ is integrable in the sense that we can 
integrate it. Actually the solutions to
the Hamilton's equation  in the variables $(I,\theta)$ are given by
$$
\left\{ \begin{array}{llll}
 \theta_{j}(t)&=&\theta_{j}(0)+ t\omega_{j},
 \quad &j\geq 1  \\
 I_{j}(t)&=&I_{j}(0), \quad &j\geq 1.
\end{array}\right.
$$
The proof of theorem \ref{birk-infini}, that we will present in section \ref{proof}, 
is very closed to the proof of 
theorem 2.6 in \cite{BDGS}. The dynamical consequences of this theorem are  
similar as those of the Birkhoff theorem in finite dimension:
\begin{corollary}\label{coro-infini}
 Assume that  $P$ belongs to the class $\T$   and that 
 $\omega$ is strongly 
 non resonant. For each $r\geq 3$ and 
 $s\geq s_{0}(r)$, there exists 
 $\e_{0}>0$ and $C>0$ such that if $\norma{(q_{0},p_{0})}=\e <\e_{0}$ the 
 solution $(q(t),p(t))$ of the Hamiltonian system associated to $H$ 
 which takes value $(q_{0},p_{0})$ at $t=0$ 
 sastisfies
\begin{itemize}
\item[(i)] $$\norma{(q_{t},p_{t})}\leq 2\e \quad \mbox{for }\va{t}\leq 
\frac{C}{\e^{r-1}}.$$
\item[(ii)] for each $j=1,\ldots,n$
$$|I_{j}(t)-I_{j}(0)|\leq \frac{\e^3}{\vaj^{2s}} \quad \mbox{for }\va{t}\leq 
\frac{C}{\e^{r-2}}$$
\item[(iii)] let $r_{1}+r_{2}=r+1$, then there exists a torus $\mathcal T_{0}\subset \Ps$ such that for each $s'< 
s-1$,
$$\mbox{dist}_{s'}((q(t),p(t),\mathcal T_{0})\leq C_{s'}\e^{r_{1}} \quad 
\mbox{for } \vat \leq 1/\e^{r_{2}}$$
where $\mbox{dist}_{s}$ denotes the distance on $\Ps$ associated with 
the norm $\norma{\cdot}_{s}$
\end{itemize}
\end{corollary}
\begin{remark}\label{kam}
This corollary remains valid for\textit{ any small initial datum} and this 
makes a big difference with the dynamical consequences of the KAM-type 
result where one has to assume that the initial datum belongs to a
Cantor-type set (cf. section \ref{KAM}). But of course the result is 
not the same, here the stability is guaranteed only for  long, \textit{but 
finite}, time. When KAM theory applies, the stability is 
inconditional, i.e. guaranteed for infinite time. Furthermore, the KAM 
theorem does not require that the perturbation be in the class $\T$.
\end{remark}
\begin{remark}
The first assertion implies in particular that the 
Hamiltonian system have \textit{almost global solutions}: if the 
initial datum is smaller than $\epsilon$ then the solution exits and 
is controled (in the initial norm) for times of order 
$\epsilon^{-r}$, the order $r$ being arbitrarily fixed at the principle. This 
consequence can be very interesting in the context of PDEs for which 
the global existence is not known (cf. \cite{BDGS}).\end{remark}

This result was first proved in \cite{BG04}. For convenience of the 
reader we repeat it here.
\proof
The proof is similar to the proof of corollary \ref{coro1},
we focus on the slight differences.
Denote $z:=(q,p)$ and $z'=\tau^{-1}(z)$ where $\tau$ is the 
transformation given by theorem \ref{birk-infini}. Then we define
$N(z):=\norma{z}_{s}^2=2\sum_{j=1}^\infty j^{2s} I_{j}(q,p)$. 
Using that $Z$ depends 
only on the actions, we have
$$
\{N,H\}(z)=\{N\circ \tau ,H\circ \tau \}\circ \tau^{-1}(z)= \{N\circ \tau, 
R\}(z').$$ 
Therefore, as in the finite dimensional case, we get
$|\dot N|\leq C N^{(r+1)/2 }$ and assertions (i) and (ii) follow.

\medskip

To prove (iii), denote by $\bar I_j:=I_j(0)$ the initial actions in the normalized
coordinates. Up to
the considered times 
\begin{equation}
\label{delI}
\left|I_j(t)-\bar I_j\right|\leq\frac{C\epsilon^{2r_1}}{j^{2s}}\ .
\end{equation}
Then, as in the proof of corollary \ref{coro2},
we define the torus
$$
\T_0:=\left\{z\in\Ps\ :\ I_j(z)=\bar I_j\ ,j\geq 1 \right\}.
$$
We have for $s'<s-1$
\begin{equation}
\label{de2}
d_{s'}(z(t),\T_0)\leq \left[\sum_{j}j^{2s'}\left| \sqrt{I_j(t)}-\sqrt{\bar
  I_j} \right|^2 \right]^{1/2}.
\end{equation}
Notice that for $a,b\geq 0$,
$$
\left|\sqrt{a}-\sqrt{b}\right|\leq \sqrt{\va{a-b}}\ .
$$ 
Thus, using
(\ref{delI}), we obtain
$$
\left[d_{s'}(z(t),\T_0) \right]^2\leq \sum_{j}\frac{j^{2s}|I_j(t)-\bar
  I_j|}{j^{2(s-s')}}\leq \sup_{j}\left(j^{2s}|I_j(t)-\bar
  I_j|\right)\sum\frac{1}{j^{2(s-s')}}
$$
which is convergent provided $s'<s-1/2$.
\endproof

\section{Application to Hamiltonian PDEs}\label{appli}
In this section we first descibe two concrete PDE's and we then verify that 
the 
abstract results of section \ref{thm} apply  to them. 

\subsection{Examples of 1-d Hamiltonian PDEs}\label{ex}
Two examples of 1-d Hamiltonian PDEs are given: the nonlinear wave 
equation and the nonlinear Schr\"odinger equation. In \cite{Cr00}, 
the reader may find much more examples 
like the Korteweg-de 
Vries equation, the Fermi-Pasta-Ulam system or the waterwaves 
system.\\
In section \ref{gene} we will  comment on recent generalisation to 
some $d$-dimensional PDE with $d\geq 2$.

\subsubsection*{Nonlinear wave equation}
As a first concrete example 
we consider a 1-d nonlinear wave equation
\begin{eqnarray}
\label{nlw}
u_{tt}-u_{xx}+V(x)u=g(x,u)\ ,\quad x\in S^1  \ , \ t\in \R \ ,
\end{eqnarray}
with Dirichlet boundary condition: $u(0,t)=u(\pi,t)=0$ for any $t$.
Here $V$ is a 2$\pi$ periodic $C^\infty$ non negative potential and $g\in
C^{\infty}(S^1\times\U)$, $\U$ being a neighbourhood of the origin in
$\R$. For compatibility reasons with the Dirichlet conditions, 
we further assume that $g(x,u)=-g(-x,-u)$ and that $V$ is even. Finally we 
assume that $g$ has a zero of order two at $u=0$ in such a way that 
$g(x,u)$ appears, in the neighborhood of $u=0$, as a perturbation 
term. 

Defining $v=u_{t}$, \eqref{nlw} reads
$$
\partial_{t}
\left(\begin{array}{c}
    u \\ v
\end{array}\right)=
\left(\begin{array}{c}
    v \\ u_{xx}-V(x)u+g(x,u)
\end{array}\right).
$$
Furthermore, let $H: H^{1}(S^1)\times L^2(S^1)\mapsto \R$ defined by
\be \label{Hnlw}
H(u,v)=\int_{S^1}\left( \frac 1 2 v^2 + \frac 1 2 u_{x}^2 +\frac 1 2 Vu^2+
G(x,u)\right) dx
\ee
where $G$ is such that
$\partial_{u}G=-g$, then \eqref{nlw} reads as an Hamiltonian system
\begin{eqnarray} \nonumber
    \partial_{t}
    \left(\begin{array}{c}
	u \\ v
    \end{array}\right)&=&
    \left(\begin{array}{cc}
	0&1 \\ -1 &0
    \end{array}\right)\left(\begin{array}{c}
    -u_{xx}+Vu+ \partial_{u}G\\ v
    \end{array}\right)\\ \label{nlwh}
    &=& J \nabla_{u,v}H(u,v)
    \end{eqnarray}
where $J=\left(\begin{array}{cc}
	0&1 \\ -1 &0
    \end{array}\right)$ represents the symplectic structure and where
    $\nabla_{u,v}=\left(\begin{array}{c}
    \nabla_{v} \\ \nabla_{v}
    \end{array}\right)$ with
    $\nabla _u$ and $\nabla _v$ denoting the $L^2$ gradient with
    respect to  $u$ and  $v$  respectively. \\
Define the operator $A:=(-\partial_{xx}+V)^{1/2}$, and introduce the
variables $(p,q)$ given by
$$
q:=A^{1/2}u\ ,\quad p:= A^{-1/2}v\ .
$$
Then, on $H^s(S^1)\times H^s(S^1)$ 
with $s\geq 1/2$, the Hamiltonian \eqref{Hnlw} takes the form $H_0+P$ with 
$$
H_0(q,p)=\frac{1}{2} \left(\left\langle Ap,p \right\rangle_{L^2} +\left\langle
Aq,q \right\rangle_{L^2} \right)
$$
and 
$$
P(q,p)=\int_{S^1}G(x,A^{-1/2}q)dx
$$
Now denote by $(\omega_j)_{j\in\Nb}$ the eigenvalues of $A$ with Dirichlet 
boundary conditions and $\phi_j$, $j\in\Nb$, the
associated eigenfunctions, i.e.
$$
A\phi_j=\omega_j\phi_j .$$
For instance, for $V=0$, we have $\phi_{j}(x)=\sin jx$ and 
$\omega_{j}= j$.\\
An element ($q,p)$ of $H^s(S^1)\times H^s(S^1)$ 
satisfying the Dirichlet 
boundary conditions  may be decomposed on the Hilbert basis 
$(\phi_{j})_{j\geq 1}$: $$q(x)=\sum_{j}q_{j}\phi_{j}(x)\quad 
\mbox{and}\quad 
p(x)=\sum_{j}p_{j}\phi_{j}(x)$$
with $(q_{j},p_{j})_{j\geq 1}\in 
\Ps=l^2_{s}(\R)\times l^2_{s}(\R)$. Then
the Hamiltonian of the non linear wave equation \eqref{nlw} reads 
on $\Ps$
$$H=\sum_{j\geq 1} \omega_{j}\frac{p_j^2+q_j^2}{2} +P$$
where $P$ belongs in $ 
C^\infty (\Ps,\R)$ and has a zero of order at least $3$ at the 
origin and $\Ps$ is endowed
with the same symplectic structure as in section \ref{model}, i.e. 
the Hamilton's equations read as in \eqref{HeqR}.

\subsubsection*{Nonlinear Schr\"odinger equation}
As a second example we consider the nonlinear Schr\"odinger equation 
\begin{eqnarray}
\label{NLS}
-\im \dot \psi=-\psi_{xx}+V\psi+\partial_{3}
  g(x,\psi,\bar \psi)\ ,\quad x\in S^1,\quad t\in \R
\end{eqnarray}
with Dirichlet boundary conditions: $\psi(0,t)=\psi(\pi,t)=0$ for any $t$. 
Here $V$ is a  $2\pi$ periodic $C^\infty$ potential. We assume that
$g(x,z_{1},z_{2})$ is $C^\infty(S^{1}\times\U)$, $\U$ being a
neighbourhood of the origin in $\C\times \C$. The notation 
$\partial_{3}$ means that we take the partial derivative with respect 
to the third argument. We also assume that $g$
has a zero of order three at $(z_{1},z_{2})=(0,0)$ and that
$g(x,z,\bar z)\in \R$.  To deal with Dirichlet boundary conditions we
have to ensure the invariance of the phase space
under the vector field associated with the equation, to this end we 
 assume that $V$ is even and that
$g(-x,-z,-\bar z) = g(x,z,\bar z)$.

Defining the Hamiltonian function of the system as 
\begin{equation}
\label{NLSH}
H=\int_{S^1}\frac{1}{2}\left(|\psi_x|^2+V|\psi|^2
\right) 
+g(x,\psi(x),\bar\psi(x))dx,
\end{equation}
equation \eqref{NLS} is equivalent to
$$
\dot \psi =i\nabla_{\bar \psi}H
$$
where $i$ represents a symplectic structure.\\
Let $A$ be the
Sturm--Liouville operator $-\partial_{xx}+V$ with Dirichlet boundary
conditions, the frequencies $\omega_j$, $j\geq 1$, are 
the corresponding
eigenvalues and the normal modes $\phi_j$  are 
the corresponding eigenfunctions. We can write $H=H_{0}+P$ with, 
for $(\psi,\bar \psi)\in H^{1}(S^1)\times H^{1}(S^1)$,
$$
H_0(\psi,\bar \psi)= \left\langle A\psi,\bar \psi \right\rangle_{L^2}
$$
and 
$$
P(\psi,\bar\psi)=\int_{S^1}g(x,\psi(x),\bar\psi(x))dx.
$$
As in the previous example an element $(\psi,\bar \psi)$ of $H^s(S^1)\times H^s(S^1)$ 
satisfying the Dirichlet 
boundary conditions  may be decomposed on the Hilbert basis 
$(\phi_{j})_{j\geq 1}$: 
$$\psi(x)=\sum_{j}\xi_{j}\phi_{j}(x) \quad \mbox{and}\quad 
\bar\psi(x)=\sum_{j}\eta_{j}\phi_{j}(x)$$
with $(\xi_{j},\eta_{j})_{j\geq 1}\in 
\Ps(\C)=l^2_{s}(\C)\times l^2_{s}(\C)$. Then
the Hamiltonian of the non linear Schr\"odinger equation \eqref{NLS} reads 
on $\Ps(\C)$
$$H=\sum_{j\geq 1} \omega_{j}\xi_{j}\eta_{j} +P.$$
Here $P$ belongs to $ 
C^\infty (\Ps,\C)$, satisfies $P(u,\bar u)\in \R$ and has a zero of order at least $3$ at the 
origin. On the other hand $\Ps(\C)$ is endowed
with the same symplectic structure as in section \ref{model}, i.e. 
the Hamilton's equations read as in \eqref{HeqC}.

Notice that defining $p$ and $q$ as the real and imaginary parts of $\psi$,
namely write $ \psi=p+\im q $ we can recover the real form \eqref{HeqR}.

\subsection{Verification of the hypothesis}
The dynamical consequences of our Birkhoff normal form theorem for PDEs are given 
in corollary \ref{coro-infini}, in particular the solution remains 
under control in the $H^s$-norm during a very long time if the 
$H^s$-norm of the initial datum is small. But this suppose that the 
Hamiltonian function of the PDE satisfies the two conditions: 
strong non resonancy of the linear frequencies and perturbation term 
in the good class. 
\subsubsection*{Verification of the condition on the perturbation term}
We work in the general framework of 1-d PDEs given in section 
\ref{model}. The Hamiltonian perturbation reads
\be \label{P}
P(q,p)=\int_{S^1}f(x,q(x),p(x))dx\ee
where $f\in C^\infty (\R^3, \R)$, $q(x)=\sum_{j\geq 1}q_{j}\phi_{j}(x)$,
$p(x)=\sum_{j\geq 1}p_{j}\phi_{j}(x)$ and\\ $((q_{j})_{j\geq 1}, 
(p_{j})_{j\geq 1})\in \Ps$. Here $(\phi_{j})_{j\geq 1}$ are the 
eigenfunctions of the selfadjoint operator $A$ and form a basis of 
the phase space.  That 
$P$ belongs to the class $\T^\nu$ is directly in relation with 
the distribution of the $\phi_{j}$'s. Actually we have
\begin{proposition} Let $\nu\geq 0$.
Assume that for each $k\geq 1$ and for each $N\geq 0$ there exists 
$C>0$ such that for all $j\in \N^k$
\be \label{phi}
\Va{\int_{S^1}\phi_{j_{1}}\ldots\phi_{j_{k}} dx}\leq C 
\frac{\mu(j)^{N+\nu}}{S(j)^N}\ee
then any $P$ of the general form \eqref{P} satisfying the symmetries 
imposed by the domain of the operator $A$ is in the class $\T^\nu$.
\end{proposition}
\proof
The Taylor's polynomial of $P$ at order $n$ reads
$$
P_{n}=\sum_{k=0}^n\sum_{(j,l)\in \Nb^{k_{1}}\times\Nb^{k_{2}} }a_{jl}q_{j_{1}}\ldots q_{j_{k_{1}}}
p_{l_{1}}\ldots p_{l_{k_{2}}}$$
with
$$
a_{jl}=\frac{1}{k_{1}! k_{2}!}\int_{S^1}\frac{\partial^k 
f}{\partial_{2}^{k_{1}}\partial_{3}^{k_{2}}}(x,0,0)
\phi_{j_{1}}(x)\ldots \phi_{j_{k_{1}}}(x)
\phi_{l_{1}}(x)\ldots \phi_{l_{k_{2}}}(x).$$
Since $P$ satisfies the symmetry conditions imposed by the domain of 
$A$, we can decompose $\frac{\partial^k 
f}{\partial_{2}^{k_{1}}\partial_{3}^{k_{2}}}(x,0,0)$ on the Hilbert 
basis $(\phi_{m})_{m\geq 1}$: 
$$\frac{\partial^k 
f}{\partial_{2}^{k_{1}}\partial_{3}^{k_{2}}}(x,0,0)=
\sum_{m\in \Nb}b_{m}\phi_{m}(x)\, .$$
Thus we get
$$
a_{jl}=\frac{1}{k_{1}! k_{2}!}\sum_{m\in \Nb}b_{m}c_{mjl}
$$
where
$$
c_{mjl}=\int_{S^1}\phi_{m}(x)
\phi_{j_{1}}(x)\ldots \phi_{j_{k_{1}}}(x)
\phi_{l_{1}}(x)\ldots \phi_{l_{k_{2}}}(x).$$
By hypothesis
$$
\va{c_{mjl}}\leq C \mu(m,j,l)^{N+\nu}S(m,j,l)^{-N}$$
where $\mu 
(m,j,l)=\mu(m,j_{1},\ldots,j_{k_{1}},l_{1},\ldots,l_{k_{2}})$ and 
$S(m,j,l)=S(m,j_{1},\ldots,j_{k_{1}},l_{1},\ldots,l_{k_{2}})$.
So it remains to verify that there exists 
$C>0$ such that 
\be \label{eq3}\sum_{m\in \Nb}|b_{m}|\mu(m,j,l)^{N+\nu}S(m,j,l)^{-N} \leq C
\mu(j,l)^{N+\nu}S(j,l)^{-N}\ .\ee
This last inequality is a consequence of the following facts:
\begin{itemize}
\item For each $i$ there exists $C_{i}$ such that $|b_{m}|(1+m)^i\leq 
C_{i}$ for all $m\geq 1$ (because $f$ is infinitely smooth and the 
$b_{m}$ act as Fourier coefficient relative to the basis $(\phi_{j})$).    
\item If $m\leq \mu(j,l)$ then $\mu(m,j,l)=\mu(j,l)$ and 
$S(m,j,l)=S(j,l)$.
\item If $m>\mu(j,l)$ then $\mu(m,j,l)<m$ and thus $\sum 
(1+m)^{-i}\mu(m,j,l)^{N+\nu}$ converges for $i$ large enough. 
\item If $m>\mu(j,l)$ then $(1+m)S(m,j,l)\geq S(j,l)$ (exercise).
\end{itemize}
\endproof

So it remains to verify condition \eqref{phi} in  concrete cases. 
We begin with a very simple one:
\begin{lemma}
If $\phi_{j}=e^{ijx}$, $j\in \Z$ then \eqref{phi} holds true with 
$\nu =0$.    
\end{lemma}
\proof
We have $\int_{S^1}\phi_{j_{1}}\ldots\phi_{j_{k}} dx= 2\pi $ if 
$j_{1}+\ldots+j_{k}=0$ and $\int_{S^1}\phi_{j_{1}}\ldots\phi_{j_{k}} 
dx=0$ if $j_{1}+\ldots+j_{k}\neq 0$. So we have to prove that there 
exists $C>0$ such that for 
any $j\in \Z^k$ satisfying $j_{1}+\ldots+j_{k}=0$,
$$
S(j)\leq C \mu(j).$$
By symmetry we can assume that $j$ is ordered, i.e. $\va{j_{1}}\leq 
\va{j_{2}}\leq
\ldots\leq\va{j_{k}}$. In this case, recall that $S(j)= 
\va{\va{j_{k}}-\va{j_{k-1}}}+\mu(j)$. But since $j_{1}+\ldots+j_{k}=0$
$$
\va{\va{j_{k}}-\va{j_{k-1}}}\leq \va{j_{k}+j_{k-1}}\leq 
\sum_{m=1}^{k-2}\va{j_{m}}\leq (k-2)\mu(j).$$
Therefore $S(j)\leq (k-1)\mu(j).$
\endproof

The condition \eqref{phi} can be verified in a much more general case:
\begin{definition}
A sequence of functions $(\phi_{j})_{j\geq 1}$ is \textbf{well localised 
with respect to the exponentials} if, writing $\phi_{j}=\sum_{l\in\Z} 
\phi_{j}^l e^{ilx}$, for each $n\geq 0$ there exists $c_{n}>0$ such 
that
\be\label{expo}
\va{\phi_{j}^l}\leq \frac{c_{n}}{\min_{\pm}(1+\va{l\pm j})^n}
\ee for all $j,l\in \Z$.
\end{definition}

\begin{example}
\label{ex.3}
If $A=-\partial_{xx}$ with Dirichlet boundary conditions,
then $\phi_{j}(x)=\sin jx$ for $j\geq 1$ which are well localized with respect
to the exponentials.
\end{example}

\begin{example}
\label{ex.5}
 Let $A=-\partial_{xx}+V$ with Dirichlet boundary conditions,
 where $V$ is a $C^{\infty}$, $2\pi$
 periodic potential.  Then $\phi_{j}(x)$ are the eigenfunctions of
 a Sturm Liouville operator. By the theory of Sturm Liouville
 operators (cf \cite{Mar86,PT}) they are well localized with respect to
 the exponentials (cf \cite{CW92}).
\end{example}
This last example applies to both PDE's we have considered in 
section \ref{ex}.
\begin{proposition}
If $(\phi_{j})_{j\geq 1}$ is well localised 
with respect to the exponentials then the condition \eqref{phi} holds 
true with $\nu =0$.    
\end{proposition}
\proof
For a 
multi-index $l\in \Z^k$, we denote $[l]=l_{1}+\ldots+l_{k}$.\\
Assume that \eqref{expo} is satisfied then
\begin{eqnarray} \nonumber \Va{\int_{S^1}\phi_{j_{1}}\ldots\phi_{j_{k}} 
dx}&=&2\pi\Va{\sum_{l\in\Z^k, 
[l]=0}\phi_{j_{1}}^{l_{1}}\ldots\phi_{j_{k}}^{l_{k}}}\\ \label{eq1}
&\leq &c_{n}^k
\sum_{l\in\Z^k, 
[l]=0}\prod_{i=1}^k \frac{1}{\min_{\pm}(1+\va{l_{i}\pm j_{i}})^n}.
\end{eqnarray} 
On the other hand, define $\e_{i}=\pm 1$ in such a way that 
$\va{l_{i}+\e_{i}j_{i}}=\min_{\pm}\va{l_{i}\pm j_{i}}$. 
By symmetry we can restrict our analysis to the case where
$j$ is ordered: $\va{j_{1}}\leq 
\ldots\leq\va{j_{k}}$. Then if $\e_{k}\e_{k-1}=-1$ we write using 
$[l]=0$,
$$
\va{j_{k}-j_{k-1}}=\va{\e_{k}j_{k}-l_{k}+\e_{k-1}j_{k-1}-l_{k-1}+
\ldots+\e_{1}j_{1}-l_{1}-\sum_{i=1}^{k-2}\e_{i}j_{i}},$$
to conclude
$$
\va{j_{k}-j_{k-1}}\leq (k-2)\mu(j)+D(l,j)
$$
where
$$
D(l,j)=\sum_{i=1}^k\va{l_{i}-\e_{i}j_{i}}.$$
If $\e_{k}\e_{k-1}=1$ we obtain similarly
$$
\va{j_{k}+j_{k-1}}\leq (k-2)\mu(j)+D(l,j).
$$
Hence, since 
$S(j)=\mu(j)+\va{j_{k}}-\va{j_{k-1}}=\mu(j)+\min_{\pm}\va{j_{k}\pm 
j_{k-1}}$, we obtain in both cases
$$
S(j)\leq (k-1)\mu(j)+D(l,j).
$$
As a consequence we have
\be\label{eq2}
\frac{\mu(j)}{S(j)}\geq \frac{1}{k-1}\ \frac{1}{1+D(l,j)}\ .
\ee
Finally notice that, by definition of $\epsilon_{i}$,
$$
\prod_{i=1}^{k}{\min_{\pm}(1+\va{l_{i}\pm j_{i}}})\geq 1+ D(l,j).
$$
Inserting this last inequality and \eqref{eq2} in \eqref{eq1} leads to
\begin{eqnarray} \nonumber
\Va{\int_{S^1}\phi_{j_{1}}\ldots\phi_{j_{k}} 
dx}&\leq&2\pi c_{n}^k
\sum_{l\in\Z^k, 
[l]=0} \frac{1}{(1+D(l,j))^n}\\ \nonumber&\leq& 
2\pi (k-1)^{N}c_{n}^k\frac{\mu(j)^N}{S(j)^N}\sum_{l\in\Z^k, 
[l]=0}\frac{1}{(1+D(l,j))^{n-N}}\ .
\end{eqnarray} 
The last sum converges for $n>N+k-1$ and thus \eqref{phi} is verified.
\endproof

\subsubsection*{Verification of the strong non resonancy condition in a 
simple case} 
This subsection is inspired by section 5 in \cite{BG04}, actually the 
case considered here is much more simple.\\ 
Let $A$ be the operator on $L^2(-\pi,\pi)$ defined by 
$$Au= -\frac{d^2u}{dx^2}+ V\star u$$
where  $V$ is a  
 $2\pi$ periodic potential and $\star$ denotes the convolution 
 product:
 $$V\star u (x) = \int_{-\pi}^{\pi}V(x-y)u(y)dy\ .$$
We consider $A$ with Dirichlet boundary conditions, i.e. 
on the domain
$D(A)$ of odd and $2\pi$-periodic $H^2$ function (cf. section 
\ref{ex}),
$$D(A)=\{u(x)= \sum_{j\geq 1}u_{j} \sin jx \mid (u_{j})_{j\geq 1}\in 
l^2_{2}(\N,\R)\}.$$
We assume that $V$ belongs to the following space ($m\geq 1$)
$$
V\in \W_{m}:=\{V(x)=\frac{1}{\pi}\sum_{j\geq 1}\frac{v_{j}}{(1+\vaj)^m}
\cos jx \mid v_{j}\in [-1/2,1/2],\ j\geq 1\}$$
that we endow with the product probability measure. Notice that a 
potential in $\W_{m}$ is in the Sobolev space $H^{m-1}$ and that we 
assume $V$ even to leave invariant $D(A)$ under the convolution 
product by $V$.

In this context the frequencies are given by 
$$\omega_{j}=j^2+\frac{v_{j}}{(1+\vaj)^m}\ ,\quad j\geq 1$$
and one has
\begin{theorem}\label{res}
    There exists a set $F_{m} \subset \W_{m}$
whose measure equals $1$ such that if $V\in F_{m}$ then the 
frequencies vector
$(\omega_{j})_{j\geq 1}$ is strongly non resonnant.
\end{theorem}
\begin{remark}\label{gen}
A similar result holds true when considering the more interesting 
case $A= -\frac{d^2}{dx^2}+ V$ which corresponds to our non linear 
Schr\"odinger equation \eqref{NLS}. But the proof is much more 
complicated (cf. \cite{BG04}).\\
In the case of the non linear wave equation \eqref{nlw} with a 
constant potential $V=m$, the frequencies reads 
$\omega_{j}=\sqrt{j^2+m}$ and it is not too difficult to prove that 
these frequencies  satisfy \eqref{A.2} for most choices of $m$ 
(see \cite{Bam03} or \cite{DS1}).
\end{remark}
Instead of proving theorem \ref{res}, we prefer to prove the 
following equivalent statement
\begin{proposition} \label{res.1}
Fix $r\geq 1$ and $\gamma >0$ small enough.  There exist positive
constants $C\equiv C_{r}$, $\alpha \equiv \alpha (r,\gamma)$, $\delta
\equiv \delta(r,\gamma)\leq \gamma$ and a set $F_{r,\gamma} \subset \W_{m}$
whose measure is larger than $1- C \gamma$ such that if $V\in F_{r,\gamma}$
then for any $N\geq 1$
\begin{eqnarray}
\label{nr.d}
\left|\sum_{j=1}^N k_{j}\omega_{j}+\epsilon_{1} \omega_{l_{1}}+\epsilon_{2}
\omega_{l_{2}}\right|
\geq \frac{\delta}{N^\alpha}
\end{eqnarray}
for any $k\in\Z^{N}$ with $\vak := \sum_{j=1}^N \va{k_{j}} \leq r$,
for any indexes 
$l_{1},\ l_{2}>N$,
and for any $\epsilon_{1},\epsilon_{2} \in \{ 0,1,-1\}$  except if\\
$k=0$ and $\epsilon_{1} =\epsilon_{2} =0$.

\end{proposition}
Theorem \ref{res} is deduced from proposition \ref{res.1} by defining 
$$
F_{m} := \cap_{r\geq 1}\cup_{\gamma >0} 
F_{r,\gamma} \ ,
$$
and remarking that this is numerable intersection of sets with full
measure. 

In order to prove proposition \ref{res.1}, we first prove that 
$\sum_{j=1}^N k_{j}\omega_{j}$ cannot accumulate on $\Z$. Precisely we have

\begin{lemma}
\label{lem.res.d}
Fix $r\geq 1$ and $\gamma >0$ small enough.
There exist positive constants $C\equiv C_{r}$, $\beta \equiv \beta (r,\gamma)$
and a set $F'_{r,\gamma} \subset \W_{m}$ whose 
measure equals $1- C \gamma$ such that if $V\in F'_{r,\gamma}$ 
then
 for any $N\geq 1$ and any $b\in \Z$
\begin{eqnarray}
\label{nr.d2}
|\sum_{j=1}^N k_{j}\omega_{j} -b|\geq  \frac{\gamma}{N^\beta}\ , 
\end{eqnarray}
for any $k\in\Z^{N}$ with $0<|k|\leq r$.
\end{lemma}

\noindent
\textit{Proof}
First notice that, given $(a_{1},\ldots,a_{r}) \neq 0$ in $\Z^r$, 
$M>0$ and 
$c\in \R$ the Lesbegue measure of 
$$
\{x\in[-M,M]^r \mid \va{\sum_{i=1}^r a_{i}x_{i} +c}<\delta \}
$$
is smaller than $(2M)^{r-1} \delta $.  Hence given $k\in \Z^N$ of
length less than $r$ and $b\in \Z$ the Lesbegue measure of
$$
\mathcal X_k:=\left\{ x\in[-1/2,1/2]^{N}\ :\ \left|\sum_{j=1}^N
k_j(j^2+x_j)+b \right|<\frac{\gamma}{N^{\beta}}  \right\}
$$
is smaller than $\gamma/N^\beta$. Now consider the set
$$
\left\{v\in[-1/2,1/2]^{N}\ :\ \left|\sum_{j=1}^N k_{j}\omega_{j} -b
\right|<\frac{\gamma}{N^\beta}\right\}, 
$$
it is contained in the set of the $v$'s such that
$(v_j/(1+|j|)^m)\in\mathcal X_k$.  The measure of this set in turn is
estimated by $(1+N)^{m(r-1)}\gamma/N^\beta$. To conclude the proof we have
to sum over all the $k$'s and the $b$'s. To count the cardinality of
the set of the $k$'s and the $b$'s to be considered remark that if
$|\sum_{j=1}^N k_{j}\omega_{j} -b|\geq \delta$ with $\delta <1$ then
$\va{b}\leq 1 + \va{\sum_{j=1}^N k_{j}\omega_{j}} \leq 1 + (1+N^{2})r$. So that  
to guarantee 
\eqref{nr.d2} for all possible choices of $k$, $b$ and $N$, it 
suffices 
to remove from $\W_{m}$ a set of measure
$$
\sum_{N\geq 1} \gamma \frac{1}{N^\beta}(1+N)^{m(r-1)}N^{r}
(1 + (1+N^{2})r) \ .
$$
Choosing $\beta := r(1+m)+4$, the last series converges and the lemma 
is proved. \endproof 

\medskip

\noindent
\textit{Proof of proposition \ref{res.1}}
First of all, for $\epsilon_{1}=\epsilon_{2}=0$, \eqref{nr.d} is a direct consequence of lemma 
\ref{lem.res.d} choosing $\alpha = \beta $, $\delta=\gamma$ and 
$F_{r,\gamma}=F'_{r,\gamma}$.\\
When $\epsilon_{1} =\pm 1$ and $\epsilon_{2}=0$, \eqref{nr.d} reads
\begin{equation}\label{l=1}
\va{\sum_{j=1}^N k_{j}\omega_{j}  \pm \omega_{l}}\geq  \frac{\delta}{N^\alpha}
\end{equation}
for some 
$l \geq N$. Notice that $\va{\sum_{j=1}^N k_{j}\omega_{j} }\leq r(N^{2}+1)$ 
and thus, if $l >2Nr$, \eqref{l=1} is always true. 
When $l\leq 2Nr$, we apply lemma \ref{lem.res.d} 
replacing $r$ by $r+1$ and $N$ by $2Nr$ to obtain \eqref{l=1} with
$\alpha =\beta (r+1,\gamma)$, $\delta = \frac{\gamma}{(2r)^\alpha}$ and  
$F_{r,\gamma}= F'_{r+1,\gamma}$. In the 
same way one proves \eqref{nr.d} when $\epsilon_{1}\epsilon_{2}=1$. So it 
remains to establish an estimate of the form
\begin{equation}\label{l=2}
\va{\sum_{j=1}^N k_{j}\omega_{j}  + \omega_{l_{1}}-\omega_{l_{2}}}
\geq \gamma \frac{\delta}{N^\alpha}
\end{equation}
for any $k\in\Z^N$, $0<|k|\leq r$ and for any $N\leq l_{1}\leq l_{2}$.\\
One has
$$
\omega_{l_{1}}-\omega_{l_{2}} = l_{1}^{2} -l_{2}^{2} 
+\frac{v_{l_{1}}}
{(1+\va{l_{1}})^m}
-\frac{v_{l_{2}}}{(1+\va{l_{2}})^m}\ .
$$
Therefore if $ 4^{\frac{1}{m}}N^{\frac{\beta +1}{m}}\gamma^{\frac{-1}{m}} 
\leq l_{1} \leq 
l_{2}$, one has with $b=l_{1}^{2} -l_{2}^{2}$
$$
\va{\omega_{l_{1}}-\omega_{l_{2}} -b} \leq   \frac{\gamma}{2N^\alpha}\ .
$$
Thus using lemma \ref{lem.res.d}, \eqref{l=2} holds true with 
$\alpha =\beta +1$, $\delta = \gamma/2$ and for $F_{r,\gamma}=F'_{r,\gamma}$.\\
Finally assume $l_{1}\leq 4^{\frac{1}{m}}N^{\frac{\beta 
+1}{m}}\gamma^{\frac{-1}{m}} $, taking into acount 
$\va{\sum_{j=1}^N k_{j}\omega_{j}}\leq r(N^{2}+1)$, \eqref{l=2} is satisfied 
when $l_{2}\geq 4^{\frac{1}{m}}N^{\frac{\beta 
+1}{m}}\gamma^{\frac{-1}{m}}3r$. So it remains to consider the case  
when  $l_{1} \leq l_{2} \leq 12rN^{\frac{\beta 
+1}{m}}\gamma^{\frac{-1}{m}}$. But in this case, we can  apply 
lemma \ref{lem.res.d} with $r$ replaced by $r+2$ and $N$ replaced by 
$12rN^{\frac{\beta 
+1}{m}}\gamma^{\frac{-1}{m}}$ to obtain \eqref{l=2} with 
$\alpha = \frac{\beta(r+2,\gamma) (\beta(r+2,\gamma)+1)}{m}$,
$\delta = \gamma(12r\gamma^{\frac{-1}{m}})^{-\beta(r+2,\gamma)}$ and 
$F_{r,\gamma}= F'_{r+2,\gamma}$.
\endproof

\section{Proof of our Birkhoff theorem in infinite 
dimension}\label{proof}
We first have to study the class of polynomials that we introduce in 
section \ref{thm}.
\subsection{Preliminary results on polynomials in 
$\T_{k}^{N,\nu}$}\label{poly}
The two propositions given in this section were first proved, in a
different context, in \cite{DS1}, \cite{DS2}. Nevertheless, for
convenience of the reader, we present slightly different proofs in our
context.

\begin{proposition}\label{tame}
    Let $k\in \Nb$,  $N\in \N$, $\nu\in[0,+\infty)$, $s\in\R$ with
    $s>\nu+3/2$,
    and let $P\in \T_{k+1}^{N,\nu}$. Then
    \begin{itemize}
    \item[(i)] $P$ extends as a continuous polynomials on $\Ps(\C)$
    and there exists a constant $C>0$ such that for all $z\in \Ps(\C)$
     $$ \va{P(z)}\leq C\norma{z}_{s}^{k+1}$$
    \item[ii)] Assume moreover that $N>s+1$, then the
    Hamiltonian vector field $X_{P}$ extends as
    a bounded function from $\Ps(\C)$ to $\Ps(\C)$.
    Furthermore, for any $s_0\in(\nu+1,s]$, there is $C>0$ such that for any
$z\in \Ps(\C)$
\be\label{eqtame}\norma{X_{P}(z)}_{s}\leq
C\norma{z}_{s}\norma{z}_{s_{0}}^{(k-1)}.\ee
    \end{itemize} 
\end{proposition}

\begin{remark}\label{rktame}
The estimate \eqref{eqtame} is of tame type (see
\cite{AlGe} for a general presentation of this concept)
and has to be 
compared with the classical tame estimate
$$
\norma{uv}_{H^s}\leq C_{s}(\norma{u}_{H^s}\norma{v}_{H^1}
+ \norma{u}_{H^1}\norma{v}_{H^s})\quad \forall u,v\in H^s(\R).$$
On the other hand, in \cite{BG04}, we obtained a Birkhoff normal form theorem
for perturbations whose Taylor's polynomials satisfies a more
general tame estimate. In this sense the theorem obtained there is
more general. 
\end{remark}

\proof  (i) Without loss of generality we can assume that $P$ is an 
homogeneous polynomial of degree $k+1$ in $\T_{k+1}^{N,\nu}$ and we write for $z\in \Ps(\C)$
\be \label{Pp}
P(z)=\sum_{j\in \Zb^{k+1}}a_{j}\, z_{j_{1}}\ldots z_{j_{k+1}}\ .
\ee
One has, using first \eqref{a} and then $\frac{\mu(j)}{S(j)}\leq 1$,
\begin{eqnarray*}
\va{P(z)}&\leq& C
\sum_{j\in 
\Zb^{k+1}}\frac{\mu(j)^{N+\nu}}{S(j)^N}\
\prod_{i=1}^{k+1}\va{z_{j_{i}}} \\
&\leq&C
\sum_{j\in 
\Zb^{k+1}}\frac{\mu(j)^\nu}{\prod_{i=1}^{k+1}\va{j_{i}}^s}\
\prod_{i=1}^{k+1}\va{j_{i}}^s\va{z_{j_{i}}} \\
&\leq&C
\sum_{j\in 
\Zb^{k+1}}\frac{1}{\prod_{i=1}^{k+1}\va{j_{i}}^{s-\nu}}\
\prod_{i=1}^{k+1}\va{j_{i}}^s\va{z_{j_{i}}} \\
&\leq&C
\left( \sum_{l\in 
\Zb}\frac{1}{\va{l}^{2s-2\nu}}\right)^{\frac{k+1}{2}}
\norma{z}_{s}^{k+1}
\end{eqnarray*}
where in the last inequality we used $k+1$ times the Cauchy-Schwarz
inequality.
Since $s>\nu+1/2$, the last sum converges and the first assertion is proved.

\medskip

\noindent (ii)  The Hamiltonian vector field of the polynomial
\eqref{Pp} reads \\
$X_{P}(z)=(\frac{\partial P}{\partial
z_{l}}(z))_{l\in \Zb}$ with, for $l\in \Zb$
$$
\frac{\partial P}{\partial z_{l}}(z)=
\sum_{j\in\Zb^{k+1}}\sum_{i=1}^{k+1}\delta_{j_{i},l}\
a_{j_{1}\ldots j_{i-1}lj_{i+1}\ldots j_{k+1}}\ z_{j_{1}}\ldots
z_{j_{i-1}}z_{j_{i+1}}\ldots z_{j_{k+1}},
$$
where $\delta_{m,n}$ denotes the Kronecker symbol.
Since the estimate \eqref{a} is symmetric with respect to
$j_{1},\ldots,j_{k+1}$ we deduce
$$
\Va{\frac{\partial P}{\partial z_{l}}(z)}\leq C\ (k+1)
\sum_{j\in 
\Zb^{k}}\frac{\mu(j,l)^{N+\nu}}{S(j,l)^N}\va{z_{j_{1}}}\ldots
\va{z_{j_{k}}}
$$
where $\mu(j,l)=\mu(j_{1},\ldots,{j_{k}},l)$ and
$S(j,l)=S(j_{1},\ldots,{j_{k}},l)$. So we obtain
\begin{eqnarray}\nonumber
\norma{X_{P}(z)}_{s}^2 &=& \sum_{l\in \Zb} l^{2s}
\Va{\frac{\partial P}{\partial z_{l}}(z)}^2 \\ \label{XP}
&\leq & C(k+1)^2\sum_{l\in \Zb} \Big(\sum_{j\in \Zb^k}
\frac{l^{s}\mu(j,l)^{N+\nu}}{S(j,l)^{N}}\va{z_{j_{1}}}\ldots
\va{z_{j_{k}}}\Big)^2 . 
\end{eqnarray}
By symmetry we may restrict ourselves to ordered multi-indices $j$,
i.e. 
satisfying $|j_{1}|\leq \ldots \leq |j_{k}|$. We then notice that
for all $l\in \Zb$ 
and for all ordered $j\in \Zb^k$
\be \label{top}
l\frac{\mu(j,l)}{S(j,l)}\leq 2 |j_{k}| .
\ee
Actually if $\val\leq 2|j_{k}|$ then \eqref{top} holds true since $
\frac{\mu(j,l)}{S(j,l)}\leq 1$. Now if $\val\geq 2|j_{k}|$
then $S(l,j)\geq |\val
-|j_{k}||\geq 1/2 \val$ and thus
$$
l\frac{\mu(j,l)}{S(j,l)}\leq 2{\mu(j,l)}\leq
2|j_{k}| $$
since $j$ is ordered. 

Fix $\epsilon >0$ such that $N-s\geq 1+\epsilon$ and $2s_{0}\geq 2\nu 
+3 +
\epsilon$. Inserting \eqref{top} in \eqref{XP} and using $\mu(j,l) \leq
S(j,l)$ we get (here $C$ designs a generic constant depending
on $k$, $N$, $\nu$, $s$ and $s_{0}$)
\begin{eqnarray*}
\norma{X_{P}(z)}_{s}^2 
&\leq & C\sum_{l\in \Zb} \Big(\sum_{j\in \Zb_{>}^k}
\frac{|j_{k}|^{s}\ \mu(j,l)^{\nu +1+\epsilon}}{S(j,l)^{1+\epsilon}}
\va{z_{j_{1}}}\ldots
\va{z_{j_{k}}}\Big)^2 .
\end{eqnarray*}
where  $\Zb_{>}^k$ denotes ther space of ordered
multi-indices.\\
Now we use that, for ordered $j$, $\mu(j,l)\leq |j_{k-1}|$ and\footnote{
for $l\geq j_{k-1}$, $S(j,l)=|j_{k-1}|+|l-j_{k}|$ and for $l\leq j_{k-1}$,
$S(j,l)\geq l$ } $S(j,l)\geq 1+|l-j_{k}|$ 
 to obtain
\begin{eqnarray*}
\norma{X_{P}(z)}_{s}^2 
&\leq & C\sum_{l\in \Zb} \Big(\sum_{j\in \Zb_{>}^k}
\frac{|j_{k}|^{s}\ |j_{k-1}|^{\nu +1+\epsilon}}{(1+|l-j_{k}|)^{1+\epsilon}}
\va{z_{j_{1}}}\ldots
\va{z_{j_{k}}}\Big)^2 \\
&=& C\sum_{l\in \Zb} \Big(\sum_{j_{k}\in \Zb} A_{j_{k}}\ 
B_{j_{k}}\Big)^2
\end{eqnarray*}
where 
\begin{eqnarray*}
A_{j_{k}}&=& |j_{k}|^s \va{z_{j_{k}}},\\
B_{j_{k}}&=& \sum_{(j_{1},\ldots,j_{k-1})\in 
\Delta_{j_{k}}} \frac{|j_{k-1}|^{\nu +1+\epsilon}}{(1+|l-j_{k}|)^{1+\epsilon}}
\prod_{i=1}^{k-1}\va{z_{j_{i}}}
\end{eqnarray*}
and $\Delta_{j_{k}}=\{ (j_{1},\ldots,j_{k-1})\in \Zb_{>}^{k-1} \mid 
j_{k-1}\leq j_{k} \}$. 
Therefore using the Cauchy-Schwarz inequality we get
$$
    \norma{X_{P}(z)}_{s}^2 \leq 
    C \norma{z}_{s}^2\sum_{l\in \Zb} \sum_{j_{k}\in \Zb}
    \frac{1}{(1+|l-j_{k}|)^{2+2\epsilon}}
\Big(    \sum_{(j_{1},\ldots,j_{k-1})\in 
\Delta_{j_{k}}}
    \prod_{i=1}^{k-1} \alpha_{j_{i}}\ \beta_{j_{i}}\Big)^2
$$
where, for $i=1,\ldots,k-1$,
$$
\alpha_{j_{i}}= |j_{i}|^{s_{0}} \va{z_{j_{i}}},\\
$$
and
\begin{eqnarray*}
\beta_{j_{i}}&=& \frac{1}{|j_{i}|^{s_{0}}},\quad \mbox{for 
}i=1,\ldots,k-2,\\
\beta_{j_{k-1}}&=& \frac{1}{|j_{k-1}|^{s_{0}-\nu-1-\epsilon}}.
\end{eqnarray*}
Then, applying $k-1$ times the Cauchy-Schwarz, we obtain  
\eqref{eqtame}.
\endproof

The second essential property of polynomials in $\T_{k}^{N,\nu}$
is captured in the following
\begin{proposition}\label{bra}
    The map $
 (P,Q)\mapsto \{P,Q\}$ define a continuous map from
 $\T_{k_{1}+1}^{N,\nu_{1}}\times\T_{k_{2}+1}^{N,\nu_{2}}$ to $
 \T_{k_{1}+k_{2}}^{N',\nu'}$ for any $N'<N-\max(\nu_{1},\nu_{2})-1$ and 
 any $\nu'>\nu_{1}+\nu_{2}+1$.
\end{proposition}
\proof
As in the proof of proposition \ref{tame}, we assume that 
$P\in \T_{k_{1}+1}^{N,\nu_{1}}$ and $Q\in \T_{k_{2}+1}^{N,\nu_{2}}$ 
are homogeneous polynomial and we write
$$
P(z)=\sum_{j\in \Zb^{k_{1}+1}}a_{j}\, z_{j_{1}}\ldots z_{j_{k_{1}+1}}
$$
 and
$$
Q(z)=\sum_{i\in \Zb^{k_{2}+1}}b_{i}\, z_{i_{1}}\ldots z_{i_{k_{2}+1}}\ .
$$
In view of the symmetry of the   estimate \eqref{a} with respect to 
the involved indices, one easily obtains
$$\{P,Q\}(z)=\sum_{(j,i)\in \Zb^{k_{1}+k_{2}}}c_{j,i}\, z_{j_{1}}\ldots 
z_{j_{k_{1}}}z_{i_{1}}\ldots z_{i_{k_{2}}}$$
whith
$$
\va{c_{j,i}}\leq (k_{1}+1)(k_{2}+1)\sum_{l\in \Zb}
\frac{\mu(j,l)^{N+\nu_{1}}}{S(j,l)^{N}}
\frac{\mu(i,l)^{N+\nu_{2}}}{S(i,l)^{N}}.
$$
Therefore it remains to prove that there exists $C>0$ such that for 
all  $j\in \Zb^{k_{1}}$ and all $i\in \Zb^{k_{2}}$,
\be \label{+}
\sum_{l\in \Zb}
 \frac{\mu(j,l)^{N+\nu_{1}}}{S(j,l)^{N}}
 \frac{\mu(i,l)^{N+\nu_{2}}}{S(i,l)^{N}}
\leq C \frac{\mu(j,i)^{N'+\nu'}}{S(j,i)^{N'}}
\ee
In order to simplify the notation, and because it does not change the estimates 
of \eqref{+},
we will assume $k_{1}=k_{2}=k$. We can also assume by symmetry that 
\begin{itemize}
\item all the indices are positive: 
$j_{1},\ldots,j_{k},i_{1},\ldots,i_{k}\geq 1$.
\item $j$ and $i$ are ordered: 
$j_{1}\leq \ldots \leq j_{k}$ and $i_{1}\leq \ldots \leq i_{k}$.
\item $j_{k}\geq i_{k}$.
\end{itemize}
Then we consider two cases: $i_{k}\geq j_{k-1}$ and 
$i_{k}\leq j_{k-1}$.

\medskip

\noindent \textbf{First case:} $\mathbf{j_{k}\geq i_{k}\geq 
j_{k-1}}$\\
We first remark that in this case,
\be \label{mu}
\mu(j,l)\leq \mu(i,j) \mbox{ and } \mu(i,l)\leq \mu(i,j).
\ee
For any multi-index $j$ we denote $\tilde 
S(j)=S(j)-\mu(j)$, i.e. $S(j)$ is the difference between the two 
largest indices. We have for all $l\in \Zb$
\be\label{S}
\St(i,j)\leq \St(i,l) +\St(j,l).\ee
Actually, if $\val \leq i_{k}$ then 
$\St(j,l)=\va{j_{k}-\max(j_{k-1},\val)}\geq j_k -i_{k}=\St(i,j)$, if 
$\val\geq j_{k}$ then $\St(i,l)=\val-i_{k}\geq j_{k}-i_{k}=\St(i,j)$ and if 
$j_{k}\leq \val\leq i_{k}$ then $\St(i,l)+\St(j,l)= 
j_{k}-i_{k}=\St(i,j)$.\\
Combining \ref{mu} and \ref{S} we get
$$
\frac{\mu(i,j)}{S(i,j)}\geq 1/2 \ \min \left(\frac{\mu(i,l)}{S(i,l)},
\frac{\mu(j,l)}{S(j,l)}\right).
$$
Assume for instance that $
\frac{\mu(j,l)}{S(j,l)}\leq\frac{\mu(i,l)}{S(i,l)}$ and let $\e>0$. We 
then have 
\begin{eqnarray*}
\sum_{l\in \Zb}\frac{\mu(j,l)^{N+\nu_{1}}}{S(j,l)^{N}}
 \frac{\mu(i,l)^{N+\nu_{2}}}{S(i,l)^{N}}&\leq& 2^{N-1-\epsilon}
 \sum_{l\in \Zb}
 {\mu(j,l)^{\nu_{1}}}
  \frac{\mu(i,j)^{N-1-\e}}{S(i,j)^{N-1-\e}}
  \frac{\mu(i,l)^{1+\e +\nu_{2}}}{S(i,l)^{1+\e}}\\
&\leq&C \frac{\mu(j,i)^{N +\nu_{1}+\nu_{2}}}{S(j,i)^{N-1-\e}} 
\end{eqnarray*}
where we used \eqref{mu} and the fact that $\sum_{l\in \Zb}
    \frac{1}{S(i,l)^{1+\epsilon}}\leq C$ uniformly with respect to $i$.
Hence in this first case, \eqref{+} holds true with $N'=N-1-\e$ and 
$\nu'= \nu_{1}+\nu_{2}+1+\e$.
    
\medskip

\noindent \textbf{Second case:} $\mathbf{j_{k}\geq j_{k-1}\geq 
i_{k}}$

\medskip

In this second case \eqref{S} remains true. Actually if
$\val\leq j_{k-1}$ then $\St(i,j)=\St(j)=\St(j,l)$ and if
$\val\geq j_{k-1}$ then
$\St(j,l)+\St(i,l)=\va{j_{k}-\val}+\va{i_{k}-\val}\geq j_{k}-j_{k-1}=\St(i,j)$.
Unfortunately \eqref{mu} may be wrong. More precisely, we still have
$\mu(i,j)\geq \mu(i,l)$ but it may occur that
$\mu(i,j)<\mu(j,l)$. Now, if we further assume that $\mu(j,l)\leq 
2\mu(i,j)$, we can proceed as in the first case to obtain the same 
estimates with an irrelevant extra factor $2^{N+\nu_{1}}$. So it 
remains to consider indices $l$ for which $\mu(j,l)> 
2\mu(i,j)$. Notice that it can occur only if $\val\geq j_{k-2}$, and 
thus we have $\mu(j,l)\leq \val$. Further, as $\mu(i,j)\geq i_{k}$, we deduce 
 $\val\geq 2i_{k}$ and thus $\va{i_{k}-l}\geq l/2$.\\
We finally have to argument differently depending wether $\St(i,l)\leq 
\St(j,l)$ or not. If it is true then, in view 
of \eqref{S}, we get
$$
S(i,j)\leq \mu(i,j)+2\St(j,l)\leq 2S(j,l).$$
Thus, using that $\mu(j,l)\leq l$, $\mu(i,l)\leq \mu(i,j)$ and
$\va{i_{k}-l}\geq l/2$,
\begin{eqnarray*}
\sum_{l,\ \St(i,l)\leq 
\St(j,l)}&&\frac{\mu(j,l)^{N+\nu_{1}}}{S(j,l)^{N}}
 \frac{\mu(i,l)^{N+\nu_{2}}}{S(i,l)^{N}}\\ &\leq& 2^{N-1-\epsilon-\nu_{1}}
 \sum_{l\in \Zb}
\frac{l^{N}}{S(i,j)^{N-1-\e-\nu_{1}}}
  \frac{\mu(i,j)^{N +\nu_{2}}}{(\va{i_{k}-l}+i_{k-1})^N}
  \frac{1}{S(j,l)^{1+\e}}\\
&\leq& C \frac{\mu(j,i)^{N +\nu_{2}}}{S(j,i)^{N-1-\e-\nu_{1}}} 
\end{eqnarray*}
where, as usual, we used  that $\sum_{l\in \Zb}
    \frac{1}{S(j,l)^{1+\epsilon}}\leq C$ uniformly with respect to 
    $j$.\\
    It remains to consider the subcase $\St(i,l)\geq 
\St(j,l)$. We then have
$$\frac{S(i,j)}{\mu(i,j)}\leq 2\frac{S(i,l)}{\mu(i,l)}$$
and thus, using again $\mu(j,l)\leq l$, $\mu(i,l)\leq \mu(i,j)$ and
$\va{i_{k}-l}\geq l/2$,
\begin{eqnarray*}
\sum_{l,\ \St(i,l)\geq 
\St(j,l)}&&\frac{\mu(j,l)^{N+\nu_{1}}}{S(j,l)^{N}}
 \frac{\mu(i,l)^{N+\nu_{2}}}{S(i,l)^{N}}\\ &\leq& 2^{N-1-\epsilon -\nu_{1}}
 \sum_{l\in \Zb}
\frac{l^{\nu_{1}}}{(\va{i_{k}-l}+i_{k-1})^{\nu_{1}}}
  \frac{\mu(i,j)^{N-1-\e-\nu_{1}}}{S(i,j)^{N-1-\e-\nu_{1}}}
  \frac{\mu(i,l)^{\nu_{1}+\nu_{2}+1+\e}}{S(i,l)^{1+\e}}\\
&\leq&C \frac{\mu(j,i)^{N +\nu_{2}}}{S(j,i)^{N-1-\e-\nu_{1}}} 
\end{eqnarray*}
Hence, in the second case, \eqref{+} holds true with $N'=N-1-\e-\nu_{1}$ and 
$\nu'= \nu_{1}+\nu_{2}+1+\e$.
\endproof
We end this section with a corollary concerning Lie transforms associated to 
polynomials in $\T^\beta$.
\begin{corollary}\label{compo}
Let $\chi$ be a real homogeneous polynomial in $\T_{l}^{\infty,\beta}$ 
with $\beta \geq 0$, $l\geq 3$  and 
denote by $\phi$ the associated Lie transform. \\
(i) Let $F\in \Hs$ with $s$ 
large enough,
 then 
$F\circ \phi\in\Hs$.\\
(ii) Let $P\in \T^{\infty,\nu}_{n} $, $\nu \geq 0$, 
$n\geq 3$ and fix $r\geq n$ an integer. Then
$$P\circ \phi = Q_{r}+ R_{r}$$
where: \\
- $Q_{r}$ is a polynomial of degree $r$ belonging to $\T^{\infty,\nu'}_{r}$ 
with\\ $\nu'= \nu +(r-n)(\beta +1)+2$,\\
- $R_{r}$ is a real Hamiltonian in the class $\T$ having a zero of order $r+1$ 
at the origin.
\end{corollary} 
\proof
(i) Let $\W$ be a neighborhood of $0$ 
in $\Ps$ such that $F$  belongs to $ 
\Hs(\W,\R)$. Since $\chi\in \T_{l}^{\infty,\beta}$, by proposition 
\ref{tame}, $\chi\in \Hs$ for $s> s_{1}=\beta +3/2$. In 
particular, for $s>s_{1}$, the flow $\Phi^t$ generated by the vector field $X_{\chi}$ 
transports an open subset of $\Ps$ into an open subset of $\Ps$.
Furthermore, since $\chi$ has a zero of order $3$, there exists $\U$ a 
neighborhood of $0$ in $\Ps$ such that the flow $\U \ni (q,p) \mapsto 
\Phi^t(q,p)\in \W$ is well defined and smooth for $0\leq t\leq 1$. By 
definition of the Lie transform, $\phi=\Phi^1$. In view of the formula 
$$X_{F\circ \phi}(q,p)=(D\phi(q,p))^{-1}X_{F}(\phi(q,p)),$$ we deduce 
that $F\circ \phi \in \Hs$ for $s>s_{1}$.\\
(ii)
We use lemma \ref{lie} (which remains valid in infinite 
dimension) to conclude
$$
\frac{d^k}{dt^k}P\circ \phi^t (q,p)\big|_{t=0}=P^{(k)}(q,p)
$$
where $P^{(k+1)}=\{P^{(k)},\chi \}$ and $P^{(0)}=P$. Therefore 
applying the Taylor's formula to $P\circ \phi^t (q,p)$ between $t=0$ 
and $t=1$ we deduce
\be\label{Pphi}
P\circ \phi(q,p)= 
\sum_{k=0}^{r-n}\frac{1}{n!}P^{(k)}(q,p)+\frac{1}{(r-n)!}\int_{0}^1 
(1-t)^r P^{(r-n+1)}(\Phi^t(q,p))dt.
\ee
Notice that $P^{(k)}(q,p)$ is a homogeneous
polynomial of degree $n+k(l-2)$ and, by 
propositions \ref{tame} and \ref{bra}, $P^{(k)}(q,p)\in \T^{\nu 
+k\beta+k+2}\cap \Hs$ for $s\geq  \nu +k\beta +k+2$. 
Therefore \eqref{Pphi} decomposes in the sum of a polynomial 
of degree $r$ in $\T^{\infty,\nu'}_{r}$ and a function in $\Hs$ having a zero 
of degree $r+1$ at the origin.

\endproof
    
\subsection{Proof of theorem \ref{birk-infini}}
We are now in position to prove theorem \ref{birk-infini}. Actually 
the proof is very close to the proof of theorem \ref{birk-fini}, i.e. 
the finite dimensional case.
So again having fixed some
$r\geq 3$, the idea is to construct iteratively for $ k=
2,\ldots,r$,  a neighborhood $\V_{k}$ of $0$ in $\Ps$ ($s$ large 
enough depending on $r$),
a canonical transformation $\tau_k$, defined on $\V_{k}$, an 
increasing sequence $(\nu_{k})_{k=2,\ldots,r}$ of positive numbers 
and real Hamiltonians $Z_{k}, P_{k+1}, Q_{k+2}, R_{k}$ such that
\begin{equation} 
\label{bb.1} 
H_{k}:=H\circ\tau_k=H_{0}+Z_{k}+P_{k+1}+Q_{k+2}+R_{k} 
\end{equation} 
and with the following properties
\bi
\item[(i)] $Z_{k}$ is a polynomial of degree $k$ in 
$\T_{k}^{\infty,\nu_{k}}$ having a zero of 
order 3 at the origin and $Z_{k}$ depends only on the (new) actions: 
$\{Z_{k},I_{j}\}=0$ for all $j\geq 1$.
\item[(ii)] $P_{k+1}$ is a homogeneous polynomial of degree $k+1$ in 
$\T_{k+1}^{\infty,\nu_{k}}$.
\item[(iii)] $Q_{k+2}$ is a polynomial of degree $r+1$ in  
$\T_{r+1}^{\infty,\nu_{k}}$ having a zero of 
order $k+2$ at the origin.
\item[(iv)] $R_{k}$ is a regular Hamiltonian belonging to 
$\Hs(\V_{k},\R)$ for $s$ large enough and having a zero of 
order $r+2$ at the origin.
\ei
First we fix $s> \nu_{r}+3/2$ to be sure to be able to apply 
proposition \ref{tame} at each step ($\nu_{r}$ will be defined later 
on independently of $s$).
Then we notice that \eqref{bb.1} at order $r$ proves theorem \ref{birk-infini} with 
$Z=Z_{r}$ and $R=P_{r+1}+R_{r}$ (notice that $Q_{r+2}=0$). In particular, 
by proposition \ref{tame} assertion (ii), $X_{P_{r+1}}$ satisfies 
\be \label{estimR}
\norma{X_{P_{r+1}}(q,p)}_{s}\leq C_{s}\norma{(q,p)}^r_{s}.
\ee
On the other hand, since $R_{r}$ 
belongs to $\Hs$, we can apply the Taylor's formula at order $r+1$ to 
$X_{R_{r}}$ to obtain the same estimate \eqref{estimR} for $R_{r}$  
on $\V\subset \V_{r}$ a neighborhood of $0$ in $\Ps$. 

\medskip

The Hamiltonian $H=H_{0}+P$ has the form 
(\ref{bb.1}) for $k=2$ with $\tau_2=I$, $\nu_{2}=\nu$, $Z_{2}=0$, 
$P_{3}$ being the Taylor's polynomial of $P$ of degree 
$3$, $Q_{4}$ being the Taylor's polynomial of $P$ of degree $r+1$ 
minus $P_{3}$ and 
$R_{2}=P-P_{3}-Q_{4}$. We 
show now how to pass from $k$ to $k+1$.

We search for $\tau_{k+1}$ of the form $\tau_{k}\circ \phi_{k+1}$, 
$\phi_{k+1}$ being a Lie transform associated to the Hamiltonian 
function $\chi_{k+1}\in \T_{k+1}^{\infty, \nu'_{k}}$ where $\nu'_{k}$ will be 
determined in lemma \ref{homolo2}. This Lie transform is well defined 
and smooth on 
a neighborhood $\V_{k+1}\subset \V_{k}$. Recall that by Taylor's formula 
we get for regular $F$
$$
F\circ \phi_{k+1}= F
+\{F,\chi_{k+1}\}+1/2 
\{\{F,\chi_{k+1}\},\chi_{k+1}\}+\ldots$$
We decompose $H_{k}\circ \phi_{k+1}$ as follows
\begin{eqnarray} 
\label{bb.6} 
H_{k}\circ \phi_{k+1}&=& H_{0}+Z_{k}+\{H_{0},\chi_{k+1}\}+P_{k+1} 
\\ 
\label{bb.7} 
&+& H_{0}\circ \phi_{k+1}-H_{0}-\{H_{0},\chi_{k+1}\}\ +\ 
Q_{k+2}\circ \phi_{k+1}
\\ 
\label{bb.8} 
&+& R_{k}\circ \phi_{k+1}\ + \ Z_{k}\circ \phi_{k+1} -Z_{k}\
+\ P_{k+1}\circ \phi_{k+1}-P_{k+1}\ .
\end{eqnarray} 
From corollary \ref{compo} and formula \eqref{Pphi}, we deduce that \eqref{bb.7}
and \eqref{bb.8} are  regular 
Hamiltonians 
having a zero of order $k+2$ at the origin and that the sum of these
terms decomposes in $P_{k+2}+Q_{k+3}+R_{k+1}$ with
$P_{k+2}$, $Q_{k+3}$ 
and $R_{k+1}$ satisfying the properties (ii), (iii) and (iv) at rank $k+1$ (with 
$\nu_{k+1}=k\nu'_{k}+ \nu_{k}+k+2$). So it remains 
to prove that $\chi_{k+1}$ can be choosen in such way that
$Z_{k+1}:=Z_{k}+\{H_{0},\chi_{k+1}\}+P_{k+1}$ 
satisfies (i). This is a consequence of the following lemma
\begin{lemma}\label{homolo2}
Let $\nu \in [0, +\infty)$ and assume that the frequencies vector of 
$H_{0}$ is strongly non resonant.
Let $Q$ be a homogeneous real  polynomial of degree $k$ in 
$\T_{k}^{\infty, \nu}$, there exist $\nu'>\nu$,  
homogeneous real polynomials $\chi$ and $Z$ of degree $k$ in $\T_{k}^{\infty, 
\nu'}$  such that
\be \label{homo2}\{H_{0},\chi\}+Q=Z\ee
and
\be \label{normal2}\{Z,I_{j}\}=0\quad \forall j\geq 1.\ee

\end{lemma}
\proof
For $j\in \Nb^{k_{1}}$ and $l\in \Nb^{k_{2}}$ with $k_{1}+k_{2}=k$ we 
denote
$$\xi^{(j)}\eta^{(l)}=\xi_{j_{1}}\ldots
\xi_{j_{k_{1}}}\eta_{l_{1}}\ldots\eta_{l_{k_{2}}}.$$
One has
$$
\{H_{0}, \xi^{(j)}\eta^{(l)}
\}
=-i\Omega(j,l) \xi^{(j)}\eta^{(l)}$$
with
$$
\Omega(j,l):= 
\omega_{j_{1}}+\ldots+\omega_{j_{k_{1}}}-\omega_{l_{1}}-\ldots 
-\omega_{l_{k_{2}}}.$$
Let $Q\in \T^{\infty, \nu}_{k}$ 
$$
Q=
\sum_{(j,l)\in \Nb^{k}}a_{jl}\xi^{(j)}\eta^{(l)}
$$
where $(j,l)\in \Nb^{k}$ means that $j\in 
\Nb^{k_{1}}$ and $l\in \Nb^{k_{2}}$ with $k_{1}+k_{2}=k$. 
Let us define
$$b_{jl}=i\Omega(j,l)^{-1}a_{ij},\quad c_{jl}=0\quad \mbox{when }
\{j_{1},\ldots,j_{k_{1}}\}\neq \{l_{1},\ldots,l_{k_{2}}\}$$
and
$$c_{jl}=a_{ij},\quad b_{jl}=0\quad \mbox{when }
\{j_{1},\ldots,j_{k_{1}}\}= \{l_{1},\ldots,l_{k_{2}}\}.$$
As $\omega$ is strongly non resonant, there exist $\gamma$ and 
$\alpha$ such that 
$$
|\Omega(j,l)|\geq \frac{\gamma}{\mu(j,l)^\alpha}
$$
for all $(j,l)\in \Nb^{k}$ with 
$\{j_{1},\ldots,j_{k_{1}}\}\neq \{l_{1},\ldots,l_{k_{2}}\}$.
Thus,
in view of definition \ref{po}, the polynomials
$$
\chi=\sum_{(j,l)\in \Nb^{k}}b_{j,l}\xi^{(j)}\eta^{(l)},
$$
and
$$
Z=\sum_{(j,l)\in \Nb^{k}}c_{j,l}\xi^{(j)}\eta^{(l)}
$$
belong in $\T_{k}^{\infty, \nu'}$ with $\nu'=\nu +\alpha$. Furthermore by 
construction they
satisfy \eqref{homo2} and \eqref{normal2}. Finally, as in the finite 
dimensional case, that $Q$ is real 
is equivalent to the symmetry relation: $\bar a_{jl}=a_{lj}$. Taking 
into acount that $\Omega_{lj}=-\Omega_{jl}$,
this symmetry remains satisfied for the polynomials 
$\chi$ and $Z$. \qed
\begin{remark}\label{rem:tame}
In this context, when we solve the so-called homological equation 
\eqref{homo2}, we loose some regularity ($\nu'= \nu +\alpha$ where 
$\alpha$ can be very large when $r$ grows). This make a big 
difference with \cite{BG04} where the tame modulus property (and a 
truncation in the Fourier modes) allowed to solve the homological 
equation in a fix space (but with growing norm).
\end{remark}

\section{Generalisations and comparison with KAM type results}
\subsection{Generalisations of theorem \ref{birk-infini}}\label{gene}
In order to apply theorem \ref{birk-infini}, the main difficulty 
consists in verifying the strong nonresonancy condition (cf definition \ref{SNR}).
When we 
consider $1$-d PDE with Dirichlet boundary conditions, this condition 
is mostly satisfied (see remark \ref{gen}). But, in a lot of other 
physical situations, 
the
condition \eqref{A.2} 
is too restrictive. Let us describe two examples:
\subsubsection*{Periodic boudary conditions}
Consider, as in section \ref{ex}, the non linear Schr\"odinger 
equation \eqref{NLS} but instead of Dirichlet 
boundary conditions, we now impose periodic boundary conditions: 
$\psi(x+2\pi,t)=\psi(x,t)$ for all $x$ and $t$ in $\R$. The 
frequencies are then the eigenvalues of
the sturm Liouville operator $A=-\partial_{xx}+ V$ with periodic 
boundary conditions. It turns out (see for instance \cite{Mar86}) that these eigenvalues can be 
indexed by $\Zb$ in such a way that $\omega_{j}= j^{2}+o(1)$, $j\in 
\Zb$. In 
particular we get $\omega_{j}-\omega_{-j}=o(1)$ and thus \eqref{A.2} cannot be
satisfied. The same problems appears with the non linear wave 
equation \eqref{nlw} with periodic boundary conditions. However 
we notice that in both cases the condition \eqref{expo} remains 
satisfied for the eigenfunctions $(\phi_{j})_{j\in \Zb}$ (see 
\cite{CW92}). That means that nonlinear terms of type \eqref{P} 
remains in the class $\T^{0}$.

\subsubsection*{Space dimension greater than 2}
Let us descibe the case of the semilinear Klein-Gordon equation on a 
sphere. Let $S^{d-1}$ be the unit sphere in $\R^d$ ($d\geq 2$) and 
$\Delta_{g}$ be the Laplace-Beltrami operator on $S^{d-1}$ for its 
canonical metric. We consider the nonlinear Klein-Gordon equation
\be \label{KG} 
(\partial_{t}^2-\Delta_{g}+m)v=-\partial_{2}f(x,v)  
\ee 
where $m$ is a strictly positive constant
and $f\in C^\infty(S^{d-1}\times \R)$ vanishes at least at 
order 3 in $v$, $\partial_{2}f$ being the derivative with respect to 
the second variable. The frequencies of the unperturbed problems are 
the square roots of the eigenvalues of the operator $-\Delta_{g}+m$:
$$\omega_{j}=\sqrt{j(j+d-2)+m}\, ,\quad j\geq 0.$$
The problem here is that 
these eigenvalues are no more simple. Denoting by $E_{j}$ the 
eigenspace associated to $\lambda_{j}=j(j+d-2)+m$, we know that 
$E_{j}$ is the space of restrictions to $S^{d-1}$ of all harmonic 
polynomials on $\R^d$ homogeneous of degree $j$. Actually for $d=2$ 
(which corresponds to 1-d nonlinear wave equation with periodic 
boundary condition), $E_{j}$ is the linear subspace spaned by 
$e^{ijx}$ and $e^{-ijx}$ and has the constant dimension two. For 
$d\geq 3$ the dimension of $E_{j}$ grows like $j^{d-1}$ (see 
\cite{BGM71} for a general reference on Laplace-Beltrami operators). Of course 
since the same frequency is now associated to different modes, 
condition \eqref{A.2} is no more satisfied. Nevertheless, if we denote 
by $e_{j}$ the dimension of $E_{j}$ and 
by $\phi_{j,l}$, $l=1,\ldots,e_{j}$ an orthonormal basis of $E_{j}$ 
then Delort and Szeftel have proved in \cite{DS1, DS2} that there 
exists $\nu \geq 0$ such that for any 
$k\geq 1$ and any $N\geq 0$ there exists a constant $C>0$ such that for any 
$j\in \N^k$ and any $l_{n}$ with $1\leq l_{n}\leq e_{j_{n}}$ ($n=1,\ldots,k$)
$$
\int_{S^{d-1}}\phi_{j_{1},l_{1}}\ldots \phi_{j_{k},l_{k}}dx\leq C 
\frac{\mu(j)^{N+\nu}}{S(j)^N}.
$$
This estimate is a  generalisation of \eqref{phi} which means 
that, generalizing the definition \ref{po}, the perturbation will 
belong to $\T^\nu$.

\subsubsection*{Generalized statement}
In this subsection we present a generalisation of theorem 
\ref{birk-infini} motivated by the previous examples. We follow the 
presentation of section \ref{sect-infini} and only focus on the new 
feature.

Fix for any $j\geq 1$ an integer $e_{j}\geq 1$. We consider the phase 
space $\Qs=\ls \times \ls$ with
$$\ls=\{(a_{j,l})_{j\geq 1,\,
1\leq l\leq e_{j}} \mid \sum_{j\geq 
1}\vaj^{2s}\sum_{l=1}^{e_{j}}\va{a_{j,l}}^2 <\infty \}
$$
that we endow with the standart norm and the standart symplectic 
structure as for $\Ps$ in section \ref{model}. We then define for 
$(q,p)\in \Qs$,
$$
H_{0}(q,p)=\frac{1}{2}\sum_{j\geq 
1}\sum_{l=1}^{e_{j}} \omega_{j,l} (q_{j,l}^2+ p_{j,l}^2)
$$
and for $j\geq 1$,
$$
J_{j}(q,p)=\frac{1}{2}\sum_{l=1}^{e_{j}}q_{j,l}^2+ p_{j,l}^2 \ .
$$
We assume that the frequencies $\omega_{j,l}$ are weakly non resonant in the 
following sense: 
\begin{definition} The vector of frequencies $(\omega_{j,l})_{j\geq 1,\,
1\leq l\leq e_{j}}$ is weakly non resonant if
for any $k\in\N$, there are $ \gamma >0$ and 
    $\alpha >0$ such that for any $j\in \Nb^{k}$, for any $l_{n}$ with
    $1\leq l_{n}\leq e_{j_{n}}$ ($n=1,\ldots,k$) and for any $1\leq i\leq 
    k$, one has 
    \begin{equation} 
    \label{A.3} 
    \left|\omega_{j_1,l_{1}}+\cdots+\omega_{j_{i},{l_{i}}}-
    \omega_{j_{i+1},l_{i+1}}-\cdots 
    -\omega_{j_{k},l_{k}} \right|\geq \frac{\gamma}{\mu (j)^{\alpha}}  
    \end{equation} 
except if $\{j_{1},\ldots,j_{i}\}=\{j_{i+1},\ldots,j_{k}\}$. 
\end{definition}
Notice that, with this definition, the frequencies of the same packet 
indexed by $j$ (i.e. $\omega_{j_1,l_{1}}$ for $1\leq l\leq e_{j}$) 
can be very close or even equal.\\
Using notations of section \ref{model}, we define
the class $\Tt_{k}^{N,\nu}$ of real polynomials of 
degree $k$ on $\Qs$
$$Q(\xi, \eta) \equiv Q(z)=\sum_{m=0}^k\sum_{j\in 
\Zb^m}\sum_{l_{1}=1}^{e_{j_{1}}}\ldots\sum_{l_{m}=1}^{e_{j_{m}}}
a_{j,l}z_{j_{1},l_{1}}\ldots z_{j_{m},l_{m}}$$
such that there exists a constant $C>0$ such that for all 
$j,l$
$$\va{a_{j,l}}\leq C \frac{\mu(j)^{N+\nu}}{S(j)^N}.$$
Then following definition \ref{T} we define a corresponding class 
$\Tt^\nu$ of Hamiltonians on $\Qs$ having a regular Hamiltonian vector 
field and Taylor's polynomials in $\Tt_{k}^{N,\nu}$.\\
Adapting the proof of theorem \ref{birk-infini} we get
\begin{theorem}\label{birk-gene}
    Assume $P\in \Tt^{\nu}$  for some $\nu \geq0$ and $\omega$ weakly 
    non resonant in the sense of \eqref{A.3}. Then 
    for any $r\geq 3$ there exists $s_{0}$ and  for any $s\geq 
    s_{0}$ there exists $\U$, $\V$ neighborhoods of the origin in 
    $\Qs$ and 
    $\tau : \V \to \U$ a real analytic canonical transformation  
     which puts $H=H_{0}+P$ in normal form up 
    to order $r$ i.e.
    $$H\circ \tau= H_{0}+Z+R$$
    with
    \bi
    \item[(i)] $Z$ is a continuous polynomial of degree $r$ which 
    commutes with all  $J_{j}$, $j\geq 1$, i.e. $\{Z, 
    J_{j}\}=0$ for all $j\geq 1$.
    \item[(ii)] $R\in C^\infty (\V,\R)$ and $\norma{X_{R}(q,p)}_{s}\leq 
    C_{s}\norma{(q,p)}_{s}^{r}$ for all $(q,p)\in \V$.
    \item[(iii)] $\tau$ is close to the identity: $\norma{\tau 
    (q,p)-(q,p)}_{s}\leq C_{s}\norma{(q,p)}_{s}^2$ for all $(q,p)\in \V$.
    \ei
\end{theorem}
This theorem is an abstract version of theorem 2.6 in \cite{BDGS}.
Notice that the concept of normal form is not the same as in theorem 
\ref{birk-infini}: the normal form $H_{0}+Z$ is no more, in general, 
integrable. The dynamical consequences are the same as in corollary 
\ref{coro-infini} but we have to replace $I_{j}$ by $J_{j}$ in the 
second assertion. Actually 
the
$J_{j}$ play the rule of almost actions: they are almost conserved 
quantities.\\ This abstract theorem applies to both examples that we present 
at the begining of this section and thus the dynamical corollary 
also. For a proof, rafinements and comments, see \cite{BG04} for the 
case of 
periodic boudary conditions and \cite{BDGS} for the case of the 
Klein-Gordon equation on the sphere. Notice that, in this last context,
the fact that $Z$ 
commutes with all the $J_{j}$ can be interpreted  saying that $H_{0}+Z$ only 
allows energy exchanges between modes in the same packet $E_{j}$ (i.e. that 
correspond to the same frequency).

We finally notice that in \cite{BDGS}, the normal form was used to prove 
an almost global 
existence result for Klein-Gordon equations with small Cauchy data on the 
sphere (and more generally on Zoll manifold). 

\subsection{Comments on KAM theory}\label{KAM}
In this section we briefly introduce the KAM theory in finite dimension 
and then we give an idea of the (partial) generalisation to the infinite dimensional 
case. Our aim is to compare these results to the Birkhoff approach 
developped in these notes.\\
For a simple introduction to the KAM theory in finite dimension we 
refer to \cite{W96} and \cite{HuIl04} (which both include a complete proof of 
KAM theorem)
and to the second chapter of \cite{KaP}. For 
infinite dimensional context, the reader may consult the books by S. 
Kuksin \cite{K1, K2} or the one by T. Kappeler and J. 
P\"oschel \cite{KaP}.

\subsubsection*{The classical KAM theorem } 
In contrast with section \ref{section:birk-fini} we consider Hamiltonian 
perturbations of \textit{Liouville} integrable system: $H=H_{0}+P$. We denote 
by\footnote{here $T^n=S^{1}\times\ldots\times S^1$, 
 $n$ times, is the $n$ dimensional torus} 
 $(I,\theta)\in \R^n\times T^n$ the action-angle variables for 
 $H_{0}$ and $\omega_{j}$, $j=1,\ldots,n$, the free frequencies. One 
 has $\omega_{j}=\frac{\partial H}{\partial I_{j}}$ and the 
 unperturbed equations read
 $$\left\{ \begin{array}{c}
 \dot I_{j}=0,\quad \ j=1,\ldots,n\ ,\\ \dot \theta_{j}= 
 \omega_{j},\quad  j=1,\ldots,n\ .
\end{array}\right.$$
The phase space $M=\U\times T^n$, where $\U$ is an open bounded domain of 
$\R^n$, is foliated by the invariant tori
$$T_{I}=T^n\times\{ I\}.$$
Our problem is to decide if these tori will persist after small 
hamiltonian perturbation 
of the system.\\
For simplicity, we assume that the perturbation is of the form 
$P=\epsilon F$. The Hamiltonian equation associated to $H$ then read
\be\label{hh}\left\{ \begin{array}{c}
\dot I_{j}=-\epsilon \frac{\partial F}{\partial \theta_{j}},\quad \ j=1,\ldots,n\ ,\\
\dot \theta_{j}= 
 \omega_{j}+\epsilon \frac{\partial F}{\partial I_{j}},\quad  j=1,\ldots,n\ .
\end{array}\right. \ee
To guarantee the persistency of $T_{I}$, it is not sufficient to assume 
the nondegenerancy of the frequencies (see definition \ref{nonresonant})
and we need the following 
\begin{definition}\label{def-dio}
A frequencies vector $\omega \in \R^n$  is diophantine 
 if there exist constants $\gamma >0$ and $\alpha>0$ such that for 
 all $0\neq k\in \Z^n$
\be\label{dio}
 \va{k\cdot \omega}\geq \frac{\gamma}{\va{k}^\alpha}\ .
\ee
\end{definition}
We denote by $D_{\gamma}$ the set of frequencies satisfying \eqref{dio} 
 for some $\alpha>0$.
It turns out that almost every vector in $\R^n$ are diophantine: by  
straighforward estimates one proves
that the set of vectors in a bouded domain of $\R^n$ 
that do not belong to $D_{\gamma}$ has Lebesgue measure $O(\gamma)$.

The second condition that we will need says that the frequencies 
effectively vary with the actions and thus we cannot stay in a 
resonnant situation when varying the actions:\\
The unperturbed system is said \textit{nondegenerate} on $\U$ if the Hessian 
matrix of $H_{0}$ 
$$\mbox{Hess}_{H_{0}}(I)=\left(\frac{\partial^{2}H_{0}}
{\partial I_{j}\partial I_{k}}(I)\right)_{1\leq 
j,k\leq n}$$ 
is invertible on $\U$. This nondegenerancy condition insures that the 
frequency map 
$$I\mapsto \omega(I)=\left(\frac{\partial H_{0}}{\partial 
I_{j}}_{1\leq j\leq n} \right)$$
is a local diffeomorphism at each point of $\U$.\\
Notice that this condition is not satisfied by the harmonic 
oscillator, $H_{0}=\sum \omega_{j} I_{j}$, for which the frequency map 
is constant. This makes difficult to directly compare theorem 
\ref{birk-fini} 
and theorem \ref{kam-fini} below.
\begin{theorem}\label{kam-fini} (The classical KAM theorem \cite{Kol, A63, 
Mos62})
Assume that $(I,\theta)\mapsto H=H_{0}+\epsilon F$ is real analytic on
the closure of $\ \U \times T^n$ and that $H_{0}$ is nondegenerate on $\U$. There exists $C>0$ 
such that if $\epsilon\leq C\gamma^2$ and if $I\in \U$ is such that
the frequencies vector $\omega(I)$
belongs to 
$D_{\gamma}$  then the corresponding torus $T_{I}$ persists  after 
perturbation.
\end{theorem}
As a dynamical consequence, we deduce that the system of equations 
\eqref{hh} has a lot of quasiperiodic solutions. But to decide if an 
invariant 
torus $T_{I}$ survives the perturbation, we have to know if the 
corresponding frequencies are in a Cantor type set. As we said in the 
introduction, this is not a realistic physical condition. That's why, 
even in the finite dimensional case, we can prefer to use the 
Birkhoff theory which provides long time stability under the condition 
that the frequencies are in an open subset of full Lebesgue measure.

The KAM theorem only concerns the Lagrangian tori, i.e. tori of 
maximal dimension. We can also wonder what happens to the lower dimensional tori. For 
instance if we fix the $n-m$ last actions to the value $0$ then we can 
define angle variables only for the $m$ first actions and the 
corresponding invariant torus $T_{I}$ is diffeomorphic to 
$T^m\times\{0_{\R^{n-m}}\}\times\{I\}$ 
whose
dimension is $m<n$. This difficult problem has been solved by H. 
Elliasson \cite{E} under the so called Melnikov condition \footnote{Actually 
V. K.Melnikov announced the result in \cite{Mel65}.} which says that, as a 
function of the first $m$ actions denoted by $\tilde I$, the quantities
$$
\sum_{j=1}^m k_{j}\omega_{j}(\tilde I) 
+\sum_{j=m+1}^nl_{j}\omega_{j}(\tilde I)
$$
does not vanish identically (and thus effectively vary with $\tilde I$ since 
$H$ is real analytic) for all non trivial $(k,l)\in \Z^m\times \Z^{n-m}$ with 
$|l|\leq 2$.\\
The theorem then says, roughly speaking, that under the hypothesis 
that $H_{0}$ is non degenerate and satisfies the Melnikov condition, 
for sufficiently small values of $\epsilon$, there exists a Cantor 
set of effective actions $\tilde I$ for which the corresponding 
invariant tori survive the small perturbation
(cf. \cite{E} or \cite{KaP} for a precise statement).

\subsubsection*{Extension to the infinite dimensional case }
When trying to extend theorem \ref{kam-fini} to the infinite 
dimensional case, we face, as in the case of Birkhoff theorem, the 
problem of extending the nonresonancy condition. It turns out that, because of 
the Dirichlet's theorem, the condition \eqref{dio} cannot be satisfied 
for all $k$
when the number of frequencies involved grows to infinity. So we 
cannot expect a polynomial control of the small divisors and it is 
very difficult to preserve tori of infinite dimension. In PDEs 
context, this would imply the existence of almost periodic solution, 
i.e. quasi-periodic solutions with a frequencies vector of infinite 
dimension. Unfortunately, up to now, there is essentially no result in this 
direction (see however the recent result by J. Bourgain \cite{Bo05b}). The only case where there exists a result applying to 
realistic PDEs concerns the perturbation of \textit{finite dimensional tori}. 
Of course, the set of finite dimensional tori is very small  within an 
infinite dimensional phase space, but it allows to describe the 
quasiperiodic solutions which is already very interesting.\\ 
A finite dimensional torus in an infinite dimensional phase space 
plays the role of a lower dimensional torus in a finite dimensional 
phase space and thus, it is not surprising that the crucial hypothesis 
in order to preserve a torus $T_{I}$ of dimension $N$
is a Melnikov condition:
\be\label{A.1}
\Va{\sum_{j=1}^N k_{j}\omega_{j} 
+\sum_{j=N+1}^{\infty}l_{j}\omega_{j}}\geq \frac{\gamma}{\vak^\alpha}\ee
for all $(k,l)\in \Z^N\times \Z^{\infty}$ with 
$|l|\leq 2$. The big differnce is that, now, the number of external 
frequencies, $\omega_{j}$ for $j\geq N+1$, is infinite. We are not 
trying to state a precise result in this direction, but it turns out 
that this Melnikov condition can be verified in certain PDE context (cf. 
\cite{K1,K2} for precise statements and further references).\\
We would like to conclude these lectures with a 
comparison of this nonresonant condition with the condition 
introduced in definition \ref{SNR}. We remark that \eqref{A.2} can be 
written in the equivalent form
\be\label{A.2ter}
\Va{\sum_{j=1}^N k_{j}\omega_{j} 
+\sum_{j=N+1}^{\infty}l_{j}\omega_{j}}\geq \frac{\gamma}{N^\alpha}
\ee
for all nontrivial $(k,l)\in \Z^N\times \Z^{\infty}$ with $\vak\leq r$ and
$|l|\leq 2$.\\
Thus, \eqref{A.1} and \eqref{A.2ter} give a control of essentially the 
same type of small divisor but, in \eqref{A.1}, $N$ (the dimension of 
the torus that we perturb) is fixed and $\vak$ (the lenght of the 
divisor that we consider) is free while, in \eqref{A.2bis}, $\vak$ 
(the degree of the monomials that we want to kill) is less than a fix 
$r$ 
and $N$ (the number of excited modes) is free.


\providecommand{\bysame}{\leavevmode\hbox to3em{\hrulefill}\thinspace}

\bigskip

{\bf Beno\^{i}t Gr\'ebert}

{\it Laboratoire de Math\'ematique Jean Leray UMR 6629,

Universit\'e de Nantes,

2, rue de la Houssini\`ere,

44322 Nantes Cedex 3, France

\smallskip
E-mail: {\tt benoit.grebert@univ-nantes.fr}
}

\end{document}